\documentclass[11pt]{amsart}

\pdfoutput=1

\usepackage{graphicx}
\usepackage[top=26mm, bottom=26mm, left=32mm, right=32mm]{geometry}
\usepackage{epsfig}
\usepackage{amsmath}
\usepackage{graphics}
\usepackage{latexsym}
\usepackage{amssymb}
\include{macros}

\def\p{\partial}
\def\xb{\boldsymbol{x}}

\def\nb{\boldsymbol{n}}
\def\nub{\boldsymbol{\tau}}

\def\({\left(}
\def\){\right)}
\def\[{\left[}
\def\]{\right]}

\numberwithin{equation}{section}
\numberwithin{figure}{section}
\theoremstyle{definition}
\newtheorem{remark}{Remark}
\numberwithin{remark}{section}

\numberwithin{definition}{section}

\newcommand{\lsp}{\vspace{3mm}}

\begin{document}

\begin{center}
\textsc{A Direct Solver for the Rapid Solution of Boundary Integral
Equations on Axisymmetric Surfaces in Three Dimensions}

\lsp

\textit{\small
Patrick Young and Per-Gunnar Martinsson\\
Dept.~of Applied Mathematics, Univ.~of Colorado at Boulder, Boulder, CO 80309-0526}

\vspace{4mm}

\begin{minipage}{130mm}\small
\noindent \textbf{Abstract:} A scheme for rapidly and accurately
computing solutions to boundary integral equations (BIEs) on rotationally
symmetric surfaces in $\mathbb{R}^{3}$ is presented. The scheme uses the
Fourier transform to reduce the original BIE defined on a surface to a
sequence of BIEs defined on a generating curve for the surface. It can
handle loads that are not necessarily rotationally symmetric. Nystr\"{o}m
discretization is used to discretize the BIEs on the generating curve.
The quadrature used is a high-order Gaussian rule that is modified near
the diagonal to retain high-order accuracy for singular kernels. The
reduction in dimensionality, along with the use of high-order accurate
quadratures, leads to small linear systems that can be inverted directly
via, \textit{e.g.}, Gaussian elimination. This makes the scheme
particularly fast in environments involving multiple right hand sides. It
is demonstrated that for BIEs associated with Laplace's equation, the
kernel in the reduced equations can be evaluated very rapidly by
exploiting recursion relations for Legendre functions. Numerical
examples illustrate the performance of the scheme; in particular, it is
demonstrated that for a BIE associated with Laplace's equation on a
surface discretized using $320\,000$ points, the set-up phase of the
algorithm takes 2 minutes on a standard desktop, and then solves can be
executed in 0.5 seconds.


\end{minipage}
\end{center}

\vspace{2mm}
\normalsize

\section{Introduction}

\begin{figure}[htbp]
	\centering
	\includegraphics[width = 0.5\linewidth]{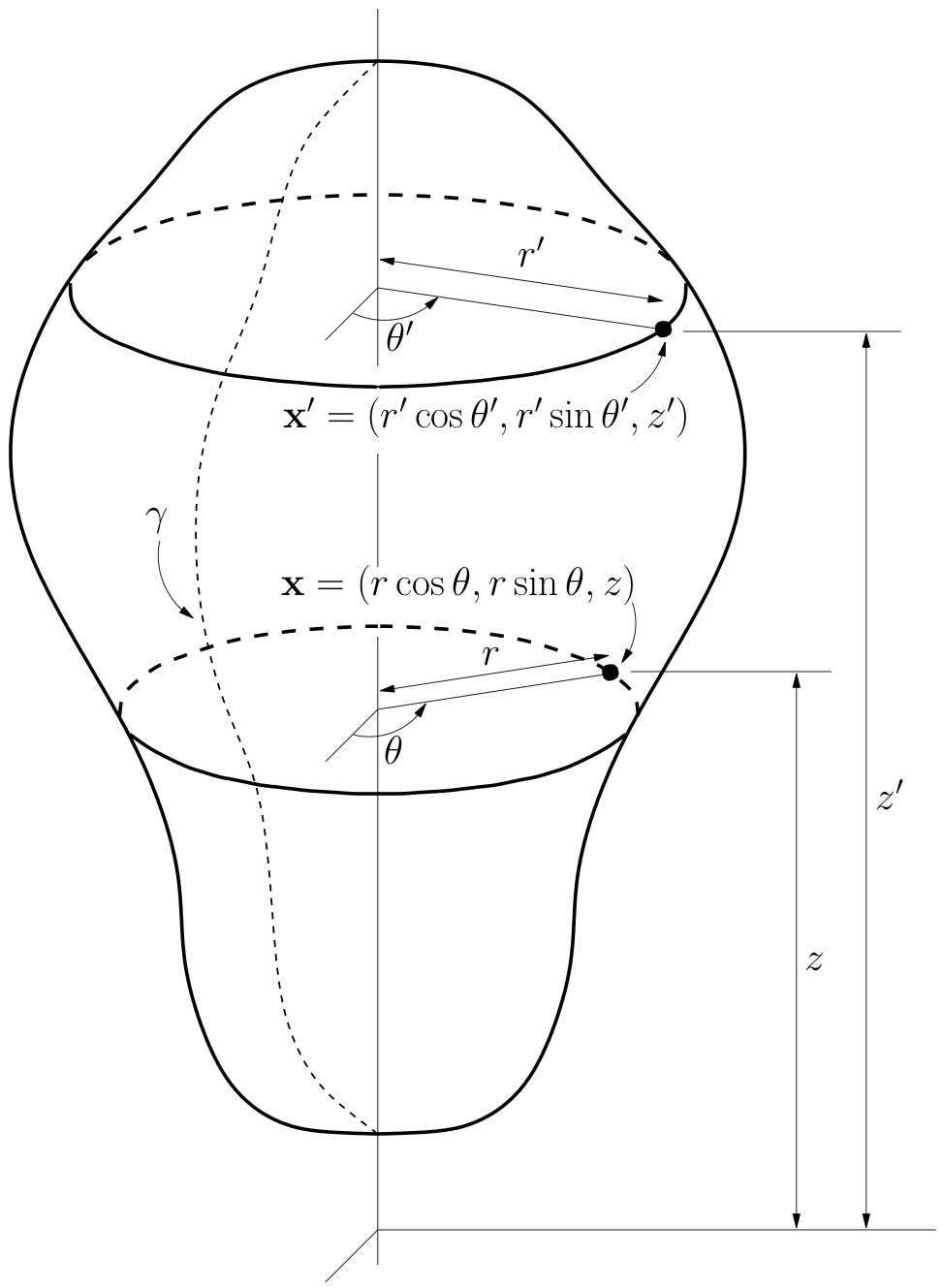}
	\caption{The axisymmetric  domain $\Gamma$ generated by the curve $\gamma$.}
	\label{fig:domain}
\end{figure}

This paper presents a numerical technique for solving boundary integral
equations (BIEs) defined on axisymmetric surfaces in $\mathbb{R}^{3}$. Specifically,
we consider second kind Fredholm equations of the form
\begin{equation}
    \sigma({\xb}) + \int_{\Gamma} k(\xb,\xb') \, \sigma(\xb') \, dA(\xb') = f(\xb), \hspace{1em} \xb \in \Gamma,
    \label{eq:fred1}
\end{equation}
under two assumptions: First, that $\Gamma$ is a surface in $\mathbb{R}^{3}$
obtained by rotating a curve $\gamma$ about an axis.
Second, that the kernel $k$ is invariant under rotation about the symmetry axis
in the sense that
\begin{equation}
    k(\xb,\xb') = k(\theta - \theta',r,z,r',z'),
    \label{eq:kerAss}
\end{equation}
where $(r,\,z,\,\theta)$ and $(r',\,z',\,\theta')$ are cylindrical
coordinates for $\xb$ and $\xb'$, respectively,
\begin{align*}
\xb  &= (r\,\cos\theta,\,r\,\sin\theta,\,z),\\
\xb' &= (r'\,\cos\theta',\,r'\,\sin\theta',\,z'),
\end{align*}
see Figure \ref{fig:domain}. Under these assumptions, the equation (\ref{eq:fred1}), which
is defined on the two-dimensional surface $\Gamma$, can via a Fourier transform in the
azimuthal variable be recast as a sequence of equations
defined on the one-dimensional curve $\gamma$. To be precise, letting $\sigma_{n}$, $f_{n}$,
and $k_{n}$ denote the Fourier coefficients of $\sigma$, $f$, and $k$, respectively (so
that (\ref{eq:forSeries1}), (\ref{eq:forSeries2}), and (\ref{eq:forSeries3}) hold), the
equation (\ref{eq:fred1}) is equivalent to the sequence of equations
\begin{equation}
\label{eq:fred2}
    \sigma_{n}(r,z) + \sqrt{2 \pi} \int_{\gamma} k_{n}(r,z,r',z') \,\sigma_{n}(r',z') \,r'\,  dl(r',z') = f_{n}(r,z),
    \hspace{1em} (r,z) \in \gamma, \hspace{0.5em} n \in \mathbb{Z}.
\end{equation}
Whenever $f$ can be represented with a moderate number of Fourier modes,
the formula (\ref{eq:fred2}) provides an efficient technique for
computing the corresponding modes of $\sigma$. The conversion of
(\ref{eq:fred1}) to (\ref{eq:fred2}) appears in, \textit{e.g.},
\cite{Rizzo:79a}, and is described in detail in Section \ref{sec:sepVar}.

Equations of the type (\ref{eq:fred1}) arise in many areas of
mathematical physics and engineering, commonly as reformulations of
elliptic partial differential equations. Advantages of a BIE
approach include a reduction in dimensionality,
often a radical improvement in the conditioning of the mathematical
equation to be solved, a natural way of handling problems defined on
exterior domains, and a relative ease in implementing high-order
discretization schemes, see, \textit{e.g.}, \cite{Atkinson:97a}.

The numerical solution of BIEs such as (\ref{eq:fred1}) poses certain
difficulties, the foremost being that the discretizations generally
involve dense matrices. Until the 1980s, this issue often times made it
prohibitively expensive to use BIE formulations as numerical tools.
However, with the advent of ``fast" algorithms (the Fast Multipole Method
\cite{rokhlin1987,rokhlin1997}, panel clustering \cite{hackbusch_1987},
etc.) for matrix-vector multiplication and the inversion of dense
matrices arising from the discretization of BIE operators, these problems
have largely been overcome for problems in two dimensions.  This is not
necessarily the case in three dimensions; issues such as surface
representation and the construction of quadrature rules in a three
dimensional environment still pose unresolved questions. The point of
recasting the single BIE (\ref{eq:fred1}) defined on a surface as the
sequence of BIEs (\ref{eq:fred2}) defined on a curve is in part to avoid
these difficulties in discretizing surfaces, and in part to exploit the
exceptionally high speed of the Fast Fourier Transform (FFT).

The reduction of (\ref{eq:fred1}) to (\ref{eq:fred2}) is only applicable
when the geometry of the boundary is axisymmetric, but presents no such
restriction in regard to the boundary load.
Formulations of this kind have been known for a long time,
and have been applied to problems in
stress analysis \cite{Bakr:85a}, scattering
\cite{Fleming:04a,Kuijpers:97a,Soenarko:93a,Tsinopoulos:99a,Wang:97a}, and potential
theory \cite{Gupta:1979a,Provatidis:98a,Rizzo:79a,Shippy:80a}.  Most of these
approaches have relied on collocation or Galerkin discretizations and
have generally relied on low-order accurate discretizations.  A complication of the axisymmetric
formulation is the need to determine the kernels $k_{n}$ for a large number
of Fourier modes $n$, since direct integration of (\ref{eq:kerAss}) through
the azimuthal variable tends to be prohibitively expensive.  When $k$
is smooth, this calculation can rapidly be accomplished using the FFT, but when $k$ is near-singular, other techniques
are required (quadrature, local refinement, etc.) that can lead to significant slowdown in
the construction of the linear systems.

The technique described in this paper improves upon previous work in
terms of both accuracy and speed. The gain in accuracy is
attained by constructing a high-order quadrature scheme for kernels
with integrable singularities. This quadrature is obtained by locally
modifying a Guassian quadrature scheme, in a manner similar to that of
\cite{Bremer:a, Bremer:b}. Numerical experiments indicate that for
simple surfaces, a relative accuracy of $10^{-10}$ is obtained
using as few as a hundred points along the generating curve.
The rapid convergence of the discretization leads to linear systems of
small size that can be solved \emph{directly} via, \textit{e.g.},
Gaussian elimination, making the algorithm particularly effective in
environments involving multiple right hand sides and when the linear
system is ill-conditioned. To describe the asymptotic complexity of
the method, we need to introduce some notation. We let $N_{\rm P}$
denote the number of panels used to discretize the generating curve
$\gamma$, we let $N_{\rm G}$ denote the number of Gaussian points in
each panel, and we let $N_{\rm F}$ denote the number of Fourier modes
included in the calculation. Splitting the computational cost into a
``set-up'' cost that needs to be incurred only once for a given geometry
and given discretization parameters, and a ``solve'' cost representing
the time required to process each right hand side, we have
\begin{equation}
\label{eq:cost_general_setup}
T_{\rm setup} \sim \underbrace{N_{\rm P}^{2}\,N_{\rm G}^{2}\,N_{\rm F}\,\log(N_{\rm F}) +
N_{\rm P}\,N_{\rm G}^{3}\,N_{\rm F}^{2}}_{\textrm{construction of linear systems}} \hspace{1.2em} +
\underbrace{N_{\rm P}^{3}\,N_{\rm G}^{3}\,N_{\rm F}}_{\textrm{inversion of systems}},
\end{equation}
and
\begin{equation}
\label{eq:cost_general_apply}
T_{\rm solve} \sim
\underbrace{N_{\rm P}\,N_{\rm G}\,N_{\rm F}\,\log(N_{\rm F})}_{\textrm{FFT of boundary data}}
\hspace{1.25em} +
\underbrace{N_{\rm P}^{2}\,N_{\rm G}^{2}\,N_{\rm F}}_{\textrm{application of inverses}}.
\end{equation}

The technique described gets particularly efficient for problems
of the form (\ref{eq:fred1}) in which the kernel $k$ is either the
single or the double layer potential associated with Laplace's equation.
We demonstrate that for such problems, it is possible to exploit recursion
relations for Legendre functions to very rapidly construct the Fourier
coefficients $k_{n}$ in (\ref{eq:fred2}). This reduces the computational
complexity of the setup (which requires the construction of a
sequence of dense matrices) from (\ref{eq:cost_general_setup}) to
$$
T_{\rm setup} \sim N_{\rm P}^{2}\,N_{\rm G}^{2}\,N_{\rm F}\,\log(N_{\rm F}) + N_{\rm P}\,N_{\rm G}^{3}\,N_{\rm F}
+ N_{\rm P}^{3}\,N_{\rm G}^{3}\,N_{\rm F}.
$$
Numerical experiments demonstrate that for a problem with $N_{\rm{P}} = 80$, $N_{\rm{G}} = 10$, and $N_{\rm F} = 400$
(for a total of $80\times 10 \times 400 = 320\,000$ degrees of freedom), this accelerated scheme
requires only 2.2 minutes for the setup, and 0.46 seconds for each solve when implemented on a
standard desktop PC.

The technique described in this paper can be accelerated further by combining
it with a fast solver applied to each of the equations in (\ref{eq:fred2}),
such as those based on the Fast Multipole Method, or the fast direct solver
of \cite{Martinsson:04a}. This would result in a highly accurate scheme
with near optimal complexity.

The paper is organized as follows:
Section \ref{sec:sepVar} describes the reduction of (\ref{eq:fred1}) to (\ref{eq:fred2})
and quantifies the error incurred by truncating the Fourier series.
Section \ref{sec:disc_BIE} presents the Nystr\"{o}m discretization of the reduced equations using high-order quadrature applicable
to kernels with integrable singularities, and the construction of the resulting linear systems.
Section \ref{sec:genAlg} summarizes the algorithm for the numerical solution of (\ref{eq:fred2}) and describes
its computational costs.
Section \ref{sec:laplace} presents the application of the algorithm for BIE formulations of Laplace's equation and describes
the rapid calculation of $k_{n}$ in this setting.
Section \ref{sec:numRes} presents numerical examples applied to problems from potential theory, and
Section \ref{sec:conclusions} gives conclusions and possible extensions and generalizations.

\section{Fourier representation of BIE}
\label{sec:sepVar}

\subsection{Problem formulation}
Suppose that $\Gamma$ is a surface in $\mathbb{R}^{3}$ obtained
by rotating a smooth contour $\gamma$ about a fixed axis and consider the
boundary integral equation
\begin{equation}
\label{eq:basic} \sigma({\xb}) + \int_{\Gamma} k(\xb,\xb') \,
\sigma(\xb') \, dA(\xb') = f(\xb), \hspace{1em} \xb \in \Gamma.
\end{equation}
In this section, we will demonstrate that if the kernel $k$ is
rotationally symmetric in a sense to be made precise, then by taking the
Fourier transform in the azimuthal variable, (\ref{eq:basic}) can be recast as a sequence of BIEs defined on the curve
$\gamma$. To this end, we introduce a Cartesian coordinate system in $\mathbb{R}^{3}$ with the
third coordinate axis being the axis of symmetry. Then cylindrical
coordinates $(r,z,\theta)$ are defined such that
\begin{align*}
x_{1} =&\ r\,\cos \theta,\\
x_{2} =&\ r\,\sin \theta,\\
x_{3} =&\ z.
\end{align*}
Figure \ref{fig:domain} illustrates the coordinate system.

The kernel $k$ in (\ref{eq:basic}) is now rotationally symmetric if
for any two points $\xb,\xb' \in \Gamma$,
\begin{equation}
\label{eq:kernel_condition}
k(\xb,\xb') = k(\theta - \theta',r,z,r',z'),
\end{equation}
where $(\theta',r',z')$ are the cylindrical coordinates of $\xb'$.


\subsection{Separation of variables}
We define for $n \in \mathbb{Z}$ the functions $f_{n}$, $\sigma_{n}$, and
$k_{n}$ via
\begin{align}
\label{eq:def_fn}
    f_{n}(r,z)            &= \int_{\mathbb{T}} \frac{e^{-in\theta}}{\sqrt{2\pi}}\, f(\theta ,r,z) \, d\theta,\\
\label{eq:def_sigman}
    \sigma_{n}(r,z)       &= \int_{\mathbb{T}} \frac{e^{-in\theta}}{\sqrt{2\pi}}\, \sigma(\theta,r,z)\, d\theta,\\
\label{eq:def_kn}
    k_{n}(r,z,r',z') &= \int_{\mathbb{T}} \frac{e^{-in\theta}}{\sqrt{2\pi}}\,k(\theta,r,z,r',z')\,d\theta.
\end{align}
The definitions (\ref{eq:def_fn}), (\ref{eq:def_sigman}), and
(\ref{eq:def_kn}) define $f_{n}$, $\sigma_{n}$, and $k_{n}$ as the
coefficients in the Fourier series of the functions $f$, $\sigma$, and
$k$ about the azimuthal variable,
\begin{align}
\label{eq:forSeries1}
f(\xb)      &= \sum_{n \in \mathbb{Z}}\frac{e^{in\theta}}{\sqrt{2\pi}}\,f_{n}(r,z),\\
\label{eq:forSeries2}
\sigma(\xb) &= \sum_{n \in \mathbb{Z}}\frac{e^{in\theta}}{\sqrt{2\pi}}\,\sigma_{n}(r,z),\\
\label{eq:forSeries3} k(\xb,\xb') = k(\theta-\theta',r,z,r',z') &=
    \sum_{n \in \mathbb{Z}}\frac{e^{in(\theta-\theta')}}{\sqrt{2\pi}}\,k_{n}(r,z,r',z').
\end{align}

To determine the Fourier representation of (\ref{eq:basic}), we multiply
the equation by $e^{-in\theta}/\sqrt{2\pi}$ and integrate $\theta$ over
$\mathbb{T}$ (for our purposes, we can think of $\mathbb{T}$ as simply
the interval $[-\pi,\,\pi]$). Equation (\ref{eq:basic}) can then be said
to be equivalent to the sequence of equations
\begin{equation}
\label{eq:step1} \sigma_{n}(r,z) + \int_{\gamma\times\mathbb{T}}
\left[\int_{\mathbb{T}}\frac{e^{-in\theta}}{\sqrt{2\pi}}\,k(\xb,\xb')\,d\theta\right]\,\sigma(\xb')\,dA(\xb')
= f_{n}(r,z),\qquad n \in \mathbb{Z}.
\end{equation}
Invoking (\ref{eq:forSeries3}), we evaluate the bracketed factor in (\ref{eq:step1}) as
\begin{multline}
\label{eq:step2}
\int_{\mathbb{T}}\frac{e^{-in\theta}}{\sqrt{2\pi}}\,k(\xb,\xb')\,d\theta
= \int_{\mathbb{T}}\frac{e^{-in\theta}}{\sqrt{2\pi}}\,k(\theta-\theta',r,z,r',z')\,d\theta\\
= e^{-in\theta'}\int_{\mathbb{T}}\frac{e^{-in(\theta-\theta')}}{\sqrt{2\pi}}\,k(\theta-\theta',r,z,r',z')\,d\theta
= e^{-in\theta'}\,k_{n}(r,z,r',z').
\end{multline}
Inserting (\ref{eq:step2}) into (\ref{eq:step1}) and executing the integration of $\theta'$ over $\mathbb{T}$, we find
that (\ref{eq:basic}) is equivalent to the sequence of equations
\begin{equation}
    \label{eq:step4}
    \sigma_{n}(r,z) + \sqrt{2\pi} \int_{\gamma}\,k_{n}(r,z,r',z')\,\sigma_{n}(r',z')\,r' \,dl(r',z') = f_{n}(r,z),
     \hspace{1em} n \in \mathbb{Z}.
\end{equation}
For future reference, we define for $n\in \mathbb{Z}$ the boundary integral operators $\mathcal{K}_{n}$ via
\begin{equation}
\label{eq:def_Kn}
[\mathcal{K}_{n} \, \sigma_{n}](r,z) = \sqrt{2\pi}\int_{\gamma} k_{n}(r,z,r',z')\,\sigma_{n}(r',z')\,r'\,dl(r',z').
\end{equation}
Then equation (\ref{eq:step4}) can be written
\begin{equation}
\label{eq:step5}
\bigl(I + \mathcal{K}_{n})\,\sigma_{n} = f_{n},\qquad n\in \mathbb{Z}.
\end{equation}
When each operator $I + \mathcal{K}_{n}$ is continuously invertible,
we can write the solution of (\ref{eq:basic}) as
\begin{equation}
\label{eq:solution_op}
\sigma(r,z,\theta)  = \sum_{n\in \mathbb{Z}} \frac{e^{in\theta}}{\sqrt{2\pi}} [(I + \mathcal{K}_{n})^{-1} f_{n}](r,z).
\end{equation}


\subsection{Truncation of the Fourier series}
When evaluating the solution operator (\ref{eq:solution_op}) in practice, we will choose a truncation parameter $N_{\rm F}$,
and evaluate only the lowest $2N_{\rm F}+1$ Fourier modes. If $N_{\rm F}$ is chosen so that
the given function $f$ is well-represented by its lowest $2N_{\rm F}+1$ Fourier modes, then in typical
environments the solution obtained
by truncating the sum (\ref{eq:solution_op}) will also be accurate. To substantiate this claim, suppose
that $\varepsilon$ is a given tolerance, and that $N_{\rm F}$ has been chosen so that
\begin{equation}
\label{eq:error_in_f}
    ||f - \sum_{n=-N_{\rm F}}^{N_{\rm F}} \frac{e^{in\theta}}{\sqrt{2\pi}} f_{n}|| \leq \varepsilon,
\end{equation}
We define an approximate solution via
\begin{equation}
\label{eq:def_sigma_approx}
\sigma_{\textrm{approx}} = \sum_{n = -N_{\rm F}}^{N_{\rm F}} \frac{e^{in\theta}}{\sqrt{2\pi}} (I + \mathcal{K}_{n})^{-1} f_{n}.
\end{equation}
From Parseval's identity, we then find that the error in the solution satisfies
\begin{align*}
|| \sigma - \sigma_{\textrm{approx}} ||^{2}
= \sum_{|n| > N_{\rm F}}||(I + \mathcal{K}_{n})^{-1} f_{n}||^{2}
\leq \sum_{|n| > N_{\rm F}}||(I + \mathcal{K}_{n})^{-1}||^{2} \, ||f_{n}||^{2} \\
\leq \left(\max_{|n| > N_{\rm F}}||(I + \mathcal{K}_{n})^{-1}||^{2}\right)\sum_{|n| > N_{\rm F}} ||f_{n}||^{2}
\leq \left(\max_{|n| > N_{\rm F}}||(I + \mathcal{K}_{n})^{-1}||^{2}\right)\varepsilon^{2}.
\end{align*}
It is typically the case that the kernel $k(\xb,\xb')$ has sufficient smoothness such that the
Fourier modes $k_{n}(r,z,r',z')$ decay as $n\rightarrow \infty$. Then
$||\mathcal{K}_{n}|| \rightarrow 0$ as $n \rightarrow \infty$ and $||(I + \mathcal{K}_{n})^{-1}|| \rightarrow 1$.
Thus, an accurate approximation of $f$ leads to an approximation in $\sigma$ that is of the same order of accuracy.
Figure \ref{fig:cond} illustrates that when $k$ is the double layer kernel associated with the Laplace
equation, and $\gamma$ is a simple curve, then $||(I + \mathcal{K}_{n})^{-1}|| \rightarrow 1$ with rapid convergence.


\section{Discretization of BIEs in two dimensions}
\label{sec:disc_BIE}

The technique in Section \ref{sec:sepVar} reduces the BIE (\ref{eq:fred1})
defined on an axisymmetric surface $\Gamma = \gamma \times \mathbb{T}$ contained in $\mathbb{R}^{3}$,
to a sequence of BIEs defined on the curve $\gamma$ contained in $\mathbb{R}^{2}$. These equations take
the form
\begin{equation}
\label{eq:2d-model}
\sigma(\xb) + \sqrt{2 \pi} \int_{\gamma} k_{n}(\xb,\xb')\,\sigma(\xb')\,r' \, dl(\xb') =
f(\xb),\qquad \xb \in \gamma,
\end{equation}
where the kernel $k_{n}$ is defined as in (\ref{eq:def_kn}).
In this section, we describe some standard
techniques for discretizing an equation such as (\ref{eq:2d-model}). For
simplicity, we limit attention to the case where $\gamma$ is a smooth
closed curve, but extensions to non-smooth curves can be handled by slight
variations of the techniques described here, \cite{Bremer:a,Bremer:b,Helsing:00a}.


\subsection{Parameterization of the curve}
\label{sec:parameterization}
Let $\gamma$ be parameterized by a vector-valued smooth function
$\nub:[0,T] \rightarrow \mathbb{R}^{2}$. The parameterization converts
(\ref{eq:2d-model}) to an integral equation defined on the interval $[0,T]$:
\begin{equation}
\sigma(\nub(t)) + \sqrt{2 \pi} \int_{0}^{T} k_{n}(\nub(t),\nub(s))\,
\sigma(\nub(s)) \, r'(\nub(s)) \, |d\nub/ds| \, ds = f(\nub(t)),
    \qquad t \in [0,T].
    \label{eq:1Dparam}
\end{equation}
To keep our formulas uncluttered, we suppress the parameterization of the
curve and the dependence on $n$ and
introduce a new kernel
\begin{equation}
\label{eq:simple_kernel}
K(t,s) = \sqrt{2 \pi} \, k_{n}(\nub(t),\nub(s))\, r'(\nub(s)) \, |d\nub/ds|,
\end{equation}
as well as the functions
\begin{equation*}
    \varphi(t) = \sigma(\nub(t)) \hspace{1em} \textrm{and} \hspace{1em} \psi(t) = f(\nub(t)).
\end{equation*}
Then techniques for solving
\begin{equation}
\label{eq:1D}
\varphi(t) + \int_{0}^{T}K(t,s)\,\varphi(s)\,ds = \psi(t),\qquad t \in
[0,T],
\end{equation}
where $\psi$ is given and $\varphi$ is to be determined, will be equally applicable to (\ref{eq:1Dparam}).


\subsection{Nystr\"{o}m method}
\label{sec:nystrom}
We will discretize (\ref{eq:1D}) via Nystr\"{om} discretization
on standard Gaussian quadrature nodes, see \cite{Atkinson:97a}.
To this end, we divide the interval $\Omega = [0,T]$ into a
disjoint partition of $N_{\rm P}$ intervals,
$$
\Omega = \bigcup_{p=1}^{N_{\rm P}} \Omega_{p},
$$
where each $\Omega_{p}$ is a subinterval called a \textit{panel}.
On each panel $\Omega_{p}$, we place the nodes of a standard $N_{\rm G}$-point Gaussian
quadrature rule $\{t_{i}^{(p)}\}_{i=1}^{N_{\rm G}}$. The idea is now to enforce (\ref{eq:1D}) at each of the $N_{\rm P}N_{\rm G}$ nodes:
$$
\sigma(t_{i}^{(p)}) + \int_{0}^{T}K(t_{i}^{(p)},s)\,\varphi(s)\,ds = \psi(t_{i}^{(p)}),
\qquad (i,p) \in \{1,2,\ldots,N_{\rm G}\} \times \{1,2,\ldots,N_{\rm P}\}.
$$
To obtain a numerical method, suppose that we can construct
for $p,q \in \{1,\,2,\,\dots,\,N_{\rm P}\}$ and $i,j \in \{1,\,2,\,\dots,\,N_{\rm G}\}$
numbers $A_{i,j}^{(p,q)}$ such that
\begin{equation}
\label{eq:Aij_function}
\int_{0}^{T}K(t_{i}^{(p)},s)\,\varphi(s)\,ds \approx
\sum_{q=1}^{N_{\rm P}}\sum_{j=1}^{N_{\rm G}} A_{i,j}^{(p,q)}\,\varphi(t_{j}^{(q)}).
\end{equation}
Then the Nystr\"{o}m method is given by solving the linear system
\begin{equation}
\label{eq:nystrom}
\varphi_{i}^{(p)} + \sum_{q=1}^{N_{\rm P}}\sum_{j=1}^{N_{\rm G}} A_{i,j}^{(p,q)}\,\varphi_{j}^{(q)} = \psi_{i}^{(p)},
\qquad (i,p) \in \{1,2,\ldots,N_{\rm G}\} \times \{1,2,\ldots,N_{\rm P}\},
\end{equation}
where $\psi_{i}^{(p)} = \psi(t_{i}^{(p)})$ and $\varphi_{i}^{(p)}$ approximates $\varphi(t_{i}^{(p)})$.
We write (\ref{eq:nystrom}) compactly as
$$
(I + A)\,\varphi = \psi
$$
where $A$ is a matrix formed by $N_{\rm P} \times N_{\rm P}$ blocks, each of size $N_{\rm G} \times N_{\rm G}$.
We let $A^{(p,q)}$ denote the block of $A$ representing the interactions between the panels $\Omega_{p}$ and $\Omega_{q}$.


\subsection{Quadrature and interpolation}
\label{sec:prelim}
We need to determine the numbers $A_{i,j}^{(p,q)}$ such that (\ref{eq:Aij_function}) holds.
The detailed construction is given in Section \ref{sec:construct_A}, and utilizes some
well-known techniques of quadrature and interpolation, which we review in this section.


\subsubsection{Standard Gaussian quadratures}
\label{sec:gauss}
Given an interval $[0,\,h]$ and a positive integer $N_{\rm G}$, the
$N_{\rm G}$-point standard Gaussian quadrature rule consists of a set of  $N_{\rm G}$ \textit{nodes}
$\{t_{j}\}_{j=1}^{N_{\rm G}} \subset [0,\,h]$, and $N_{\rm G}$ \textit{weights} $\{w_{j}\}_{j=1}^{N_{\rm G}}$ such that
$$
\int_{0}^{h} g(s)\,ds = \sum_{j=1}^{N_{\rm G}}w_{j}\,g(t_{j}),
$$
whenever $g$ is a polynomial of degree at most $2N_{\rm G}-1$, and such that
$$
\int_{0}^{h} g(s)\,ds = \sum_{j=1}^{N_{\rm G}}w_{j}\,g(t_{j}) + O(h^{2N_{\rm G}}),
$$
whenever $g$ is a function with $2N_{\rm G}$ continuous derivatives, see \cite{Abramowitz:65a}.


\subsubsection{Quadrature rules for singular functions}
\label{sec:modGauss}
Now suppose that given an interval $[0,\,h]$ and a point $t \in [-h,\,2h]$, we seek to integrate
over $[0,\,h]$ functions $g$ that take the form
\begin{equation}
\label{eq:kForm}
g(s) = \phi_{1}(s) \log |s - t| + \phi_{2}(s),
\end{equation}
where $\phi_{1}$ and $\phi_{2}$ are polynomials of degree at most $2N_{\rm G}-1$.
Standard Gaussian quadrature would be highly inaccurate if applied to
integrate (\ref{eq:kForm}).  Rather, seek a $N_{\rm{G}}'$-node quadrature
that will evaluate
\begin{equation}
\label{eq:logInt}
\int_{0}^{h} g(s) \, ds
\end{equation}
exactly. Techniques for constructing such generalized quadratures are
readily available in the literature, see for example \cite{Kolm:01a}.
These quadratures will be of degree $2N_{\rm G}-1$, just as with standard
Gaussian quadratures and exhibit comparable accuracy, although in general
$N_{\rm{G}}' > N_{\rm G}$. The generalized quadratures used in this paper were
determined using the techniques of \cite{Kolm:01a}, and can be found in the appendix.

We observe that the quadrature nodes constructed by such methods are
typically different from the nodes of the standard Gaussian quadrature.
This complicates the the construction of the matrix $A$, as described in
Section \ref{sec:construct_A}.


\subsubsection{Lagrange interpolation}
\label{sec:lagrange}
Let $\{t_{j}\}_{j=1}^{N_{\rm G}}$ denote the nodes of a $N_{\rm G}$-point Gaussian quadrature
rule on the interval $[0,h]$. If the values of a polynomial $g$ of degree at most $N_{\rm G}-1$
are specified at these nodes, the entire polynomial $g$ can be recovered via
the formula
$$
g(s) = \sum_{j=1}^{N_{\rm G}} L_{j}(s)\,g(t_{j}),
$$
where the functions $L_{j}$ are the Lagrange interpolating polynomials
$$
L_{j}(s) = \prod_{i \neq j} \( \frac{s - t_{i}}{t_{j} - t_{i}} \).
$$
If $g$ is a smooth function with $N_{\rm G}$ continuous derivatives that is not a polynomial, then the Lagrange interpolant
provides an approximation to $g$ satisfying
$$
\left|g(s) - \sum_{j=1}^{N_{\rm G}} L_{j}(s)\,g(t_{j})\right| \leq C\,h^{N_{\rm G}},
$$
where
$$
C = \( {\displaystyle\sup_{s \in [0,h]}} |g^{(N_{\rm G})}(s)| \) / N_{\rm G}!.
$$


\subsection{Constructing the matrix $A$}
\label{sec:construct_A}
Using the tools reviewed in Section \ref{sec:prelim}, we are now in position to construct
numbers $A_{i,j}^{(p,q)}$ such that (\ref{eq:Aij_function}) holds. We first note
that in forming block $A^{(p,q)}$ of $A$, we need to find numbers $A_{i,j}^{(p,q)}$ such that
\begin{equation}
\label{eq:logitech0}
\int_{\Omega_{q}} K(t_{i}^{(p)},s)\,\varphi(s)\,ds \approx \sum_{j=1}^{N_{\rm G}} A_{i,j}^{(p,q)}\,\varphi(t_{j}^{(q)}),
\qquad i = 1,\,2,\,\dots,\,N_{\rm G}.
\end{equation}
When $\Omega_{p}$ and $\Omega_{q}$ are well separated, the integrand in (\ref{eq:logitech0})
is smooth, and our task is easily solved using standard Gaussian quadrature
(as described in Section \ref{sec:gauss}):
$$
\int_{\Omega_{q}}K(t_{i}^{(p)},s)\,\varphi(s)\,ds
\approx \sum_{j=1}^{N_{\rm G}} w_{j}\,K(t_{i}^{(p)},t_{j}^{(q)})\,\varphi(t_{j}^{(q)}).
$$
It directly follows that the $ij$ entry of the block $A^{(p,q)}$ takes the form
\begin{equation}
\label{eq:formula_offd}
A^{(p,q)}_{i,j} = w_{j}\,K(t_{i}^{(p)},\,t_{j}^{(q)}).
\end{equation}

Complications arise when we seek to form a diagonal block $A^{(p,p)}$, or even
a block that is adjacent to a diagonal block. The difficulty is that the
kernel $K(t,s)$ has a singularity as $s\rightarrow t$.
To be precise, for any fixed $t$, there exist smooth
functions $u_{t}$ and $v_{t}$ such that
$$
K(t,s) = \log|t-s|\,u_{t}(s) + v_{t}(s).
$$
We see that when $t_{i}^{(p)}$ is a point in $\Omega_{q}$ the integrand in (\ref{eq:logitech0}) becomes singular.
When $t_{i}^{(p)}$ is a point in a panel neighboring $\Omega_{q}$, the problem is
less severe, but Gaussian quadrature would still be inaccurate.
To maintain full accuracy, we use the modified quadrature rules
described in Section \ref{sec:modGauss}. For every node $t_{i}^{(p)} \in \Omega_{p}$,
we construct a quadrature $\{\hat{w}_{i,\ell}^{(p,q)},\hat{t}^{(p,q)}_{i,\ell}\}_{\ell=1}^{N_{\rm{G}}'}$
such that
\begin{equation}
\label{eq:logitech1}
\int_{\Omega_{q}}K(t_{i}^{(p)},s)\,\varphi(s)\,ds
\approx \sum_{\ell=1}^{N_{\rm{G}}'} \hat{w}_{i,\ell}^{(p,q)}\,K(t_{i}^{(p)},\hat{t}_{i,\ell}^{(p,q)})\,\varphi(\hat{t}_{i,\ell}^{(p,q)}).
\end{equation}
In order to have a quadrature evaluated at the Gaussian nodes $t_{j}^{(q)} \in \Omega_{q}$,
we next use Lagrange interpolation as described in Section \ref{sec:lagrange}. With $\{L_{j}^{(q)}\}_{j=1}^{N_{\rm G}}$
denoting the Lagrange interpolants of order $N_{\rm G}-1$ defined on $\Omega_{q}$, we have
\begin{equation}
\label{eq:logitech2}
\varphi(t) \approx \sum_{j=1}^{N_{\rm G}} L_{j}^{(q)}(t)\,\varphi(t_{j}^{(q)}),\qquad t \in \Omega_{q}.
\end{equation}
Inserting (\ref{eq:logitech2}) into (\ref{eq:logitech1}), we find that
$$
\int_{\Omega_{q}}K(t_{i}^{(p)},s)\,\varphi(s)\,ds \approx
\sum_{\ell=1}^{N_{\rm{G}}'}\hat{w}_{i,\ell}^{(p,q)}\,K(t_{i}^{(p)},\hat{t}_{i,\ell}^{(p,q)})\,\sum_{j=1}^{N_{\rm G}}L_{j}^{(q)}(\hat{t}_{i,\ell}^{(p,q)})\,
\varphi(t_{j}^{(q)}).
$$
We now find that the block $A^{(p,q)}$ of $A$ has entries
\begin{equation}
\label{eq:logitech3}
A^{(p,q)}_{i,j} = \sum_{\ell=1}^{N_{\rm{G}}'}
\hat{w}_{i,\ell}^{(p,q)}\,K(t_{i}^{(p)},\hat{t}_{i,\ell}^{(p,q)})\,L_{j}^{(q)}(\hat{t}_{i,\ell}^{(p,q)}),
\qquad i,j \in \{1,2,\ldots,N_{G} \}.
\end{equation}
We observe that the formula (\ref{eq:logitech3}) is quite expensive to evaluate; in addition to the summation,
it requires the construction of a quadrature rule for each point $t_{i}^{(p)}$ and evaluation of Lagrange interpolants.
Fortunately, this process must be executed for at most three blocks in each row of blocks of $A$.


\section{A general algorithm}
\label{sec:genAlg}

\subsection{Summary}
\label{sec:summary}
At this point, we have shown how to convert a BIE defined on
an axisymmetric surface in $\mathbb{R}^{3}$ to a sequence of equations defined on
a curve in $\mathbb{R}^{2}$ (Section \ref{sec:sepVar}), and then how to discretize
each of these reduced equations (Section \ref{sec:disc_BIE}). Putting everything
together, we obtain the following algorithm for solving (\ref{eq:fred1}):

\lsp

\begin{enumerate}
\item Given the right hand side $f$, and a computational tolerance $\varepsilon$,
determine a truncation parameter $N_{\rm F}$ such that (\ref{eq:error_in_f}) holds.

\lsp

\item Form for $n = -N_{\rm F},\,-N_{\rm F} \,+1,\,-N_{\rm F}+2,\,\dots,\,N_{\rm F}$ the matrix $A_{n}$ discretizing
the equation (\ref{eq:step5}) encapsulating the $n$'th Fourier mode. The matrix is
formed via Nystr\"{o}m discretization as described in Section \ref{sec:disc_BIE} with
the discretization parameters $N_{\rm P}$ and $N_{\rm G}$ chosen to meet the computational tolerance
$\varepsilon$.

\lsp

\item
Evaluate via the FFT the terms $\{f_{n}\}_{n=-N_{\rm F}}^{N_{\rm F}}$  in the Fourier representation of $f$
(as defined by (\ref{eq:def_fn})),
and solve for $n = -N_{\rm F},\,-N_{\rm F}+1,\,-N_{\rm F}+2,\,\dots,\,N_{\rm F}$  the equation
$(I + A_{n})\,\sigma_{n} = f_{n}$ for $\sigma_{n}$. Construct $\sigma_{\rm approx}$
using formula (\ref{eq:def_sigma_approx}) evaluated via the FFT.

\end{enumerate}

\lsp

The construction of the matrices $A_{n}$ in Step 2 can be accelerated using the FFT
(as described in Section \ref{sec:fastform_An}), but even with such acceleration, it
is typically by a wide margin the most expensive part of the algorithm. However, this
step needs to be performed only once for any given geometry,
and given discretization parameters $N_{\rm F}$, $N_{\rm P}$, and $N_{\rm G}$.
The method therefore becomes particularly efficient when (\ref{eq:fred1}) needs to be
solved for a sequence of right-hand sides.
In this case, it may be worth the cost to pre-compute the inverse of each matrix $I+A_{n}$.


\subsection{Techniques for forming the matrices}
\label{sec:fastform_An}
We need to construct for each Fourier mode $n$, a matrix $A_{n}$ consisting
of $N_{\rm P} \times N_{\rm P}$ blocks $A_{n}^{(p,q)}$, each of size $N_{\rm G} \times N_{\rm G}$. Constructing an off-diagonal block $A_{n}^{(p,q)}$ when $\Omega_{p}$ and $\Omega_{q}$ are not directly adjacent is relatively
straightforward. For any pair of nodes $t_{i}^{(p)} \in \Omega_{p}$ and $t_{j}^{(q)} \in \Omega_{q}$,
we need to construct the numbers, \textit{cf.}~(\ref{eq:simple_kernel}) and (\ref{eq:formula_offd}),
\begin{equation}
\label{eq:A_cong_offd}
A^{(p,q)}_{n;i,j} = \sqrt{2 \pi} \, w_{j}\,k_{n}(\nub(t_{i}^{(p)}),\,\nub(t_{j}^{(q)}))\,r'(\nub(t_{j}^{(q)})) \, |d\nub(t_{j}^{(q)})/ds|,
\end{equation}
for $n = -N_{\rm F},\,-N_{\rm F}+1,\,\dots,\,N_{\rm F}$, where $\nub$ is a parameterization of $\gamma$
(see Section \ref{sec:parameterization}) and
the kernel $k_{n}$ is defined by (\ref{eq:def_kn}). Fortunately, we do not need to
explicitly evaluate the integrals in (\ref{eq:def_kn}) since all the $2N_{\rm F}+1$ numbers can
be evaluated by a single application of the FFT to the function
\begin{equation}
\label{eq:theta_function}
\theta \mapsto k(\theta,\,\nub(t_{i}^{(p)}),\,\nub(t_{j}^{(q)})).
\end{equation}
When $\nub(t_{i}^{(p)})$ is not close to $\nub(t_{j}^{(q)})$, the function in (\ref{eq:theta_function})
is smooth, and the trapezoidal rule implicit in applying the FFT is highly accurate.

Evaluating the blocks on the diagonal, or directly adjacent to the diagonal is
somewhat more involved. The matrix entries are now given by the formula,
\textit{cf.}~(\ref{eq:simple_kernel}) and (\ref{eq:logitech3}),
\begin{equation}
A^{(p,q)}_{k;i,j} = \sum_{\ell=1}^{N_{\rm{G}}'}
\hat{w}_{i,\ell}^{(p,q)}\,k_{n}(\nub(t_{i}^{(p)}),\nub(\hat{t}_{i,\ell}^{(p,q)}))\,
r'(\hat{t}_{i,\ell}^{(p,q)}) \, |d\nub(\hat{t}_{i,\ell}^{(p,q)})/ds|\,
L_{j}^{(q)}(\hat{t}_{i,\ell}^{(p,q)}),
\end{equation}
where $\nub$ and $k_{n}$ are as in (\ref{eq:A_cong_offd}). To further complicate things,
the points $\nub(t_{i}^{(p)})$ and $\nub(\hat{t}_{i,\ell}^{(p,q)})$ are now in close proximity to
each other, and so the functions
\begin{equation}
\label{eq:theta_function2}
\theta \mapsto k(\theta,\,\nub(t_{i}^{(p)}),\,\nub(\hat{t}_{i,\ell}^{(p,q)}))
\end{equation}
have a sharp peak around the point $\theta = 0$. They are typically still easy to integrate
away from the origin, so the integrals in (\ref{eq:def_kn}) can for a general kernel be
evaluated relatively efficiently using quadratures that are adaptively refined near
the origin.

Even with the accelerations described in this section, the cost of forming the matrices
$A_{n}$ tends to dominate the computation whenever the kernels $k_{n}$ must be evaluated
via formula (\ref{eq:def_kn}). In particular environments, it is possible to side-step
this problem by evaluating the integral in (\ref{eq:def_kn}) analytically. That this can
be done for the single and double layer kernels associated with Laplace's equation is
demonstrated in Section \ref{sec:laplace}.


\subsection{Computational costs}
The asymptotic cost of the algorithm described in Section \ref{sec:summary}
has three components: (a) the cost of forming the matrices $\{A_{n}\}_{n=-N_{\rm F}}^{N_{\rm F}}$,
(b) the cost of transforming functions from physical space to Fourier space and back,
and (c) the cost of solving the linear systems $(I + A_{n})\,\sigma_{n} = f_{n}$.
In this section, we investigate the asymptotic cost of these steps. We consider
a situation where $N_{\rm F}$ Fourier modes need to be resolved, and where
$N_{\rm P}\times N_{\rm G}$ nodes are used to discretize the curve $\gamma$.

\lsp

\noindent\textit{(a) Cost of forming the linear systems:} Suppose first that we
have an analytic formula for each kernel $k_{n}$. (As we do, \textit{e.g.}, when
the original BIE (\ref{eq:fred1}) involves either the single or the double layer kernel
associated with Laplace's equation, see Section \ref{sec:laplace}.) Then the cost $T_{\rm mat}$
of forming the matrices satisfies
$$
T_{\rm mat} \hspace{2em} \sim
\underbrace{N_{\rm P}^{2}\,N_{\rm G}^{2}\,N_{\rm F}}_{\textrm{cost from kernel evaluations}} \hspace{0.5em} +
\underbrace{N_{\rm P}\,N_{\rm G}^{3}\,N_{\rm F}}_{\textrm{cost from composite quadrature}}.
$$
When the kernels have to be evaluated numerically via formula (\ref{eq:def_kn}),
the cost of forming the matrices is still moderate.  In the rare situations where
the kernel is smooth, standard Gaussian quadrature can be used everywhere
and the FFT acceleration described in Section \ref{sec:fastform_An}
can be used for all entries. In this situation,
$$
T_{\rm mat} \sim N_{\rm P}^{2}\,N_{\rm G}^{2}\,N_{\rm F}\,\log(N_{\rm F}).
$$
In the more typical situation where each kernel $k_{n}$ involves an integrable
singularity at the diagonal, the FFT acceleration can still be used to rapidly
evaluate all entries well-removed from the diagonal. However, entries close to
the diagonal must be formed via the composite quadrature rule combined with
numerical evaluation of $k_{n}$ via an adaptive quadrature. In this situation,
$$
T_{\rm mat} \sim
N_{\rm P}^{2}\,N_{\rm G}^{2}\,N_{\rm F}\,\log(N_{\rm F}) +
N_{\rm P}\,N_{\rm G}^{3}\,N_{\rm F}^{2}.
$$

\lsp

\noindent\textit{(b) Cost of Fourier transforms:} The boundary data defined on the surface
must be converted into the Fourier domain.  This is executed via the FFT
at a cost $T_{\rm fft}$ satisfying
\begin{equation}
\label{eq:T_fft}
T_{\rm fft} \sim N_{\rm P}\,N_{\rm G}\,N_{\rm F}\,\log(N_{\rm F}).
\end{equation}
We observe that the constant of proportionality in (\ref{eq:T_fft}) is very small,
and the cost of this step is typically negligible compared to the costs of the
other steps.

\lsp

\noindent\textit{(c) Cost of linear solves:} Using standard Gaussian elimination,
the cost $T_{\rm solve}$ of solving $N_{\rm F}$ linear systems $(I + A_{n})\,\sigma_{n} = f_{n}$,
each of size $N_{\rm P}\,N_{\rm G} \times N_{\rm P}\,N_{\rm G}$, satisfies
$$
T_{\rm solve} \sim N_{\rm P}^{3}\,N_{\rm G}^{3}\,N_{\rm F}.
$$
In situations where the equations need to be solved for multiple right hand sides,
it pays off to first compute the inverses $(I + A_{n})^{-1}$, and then simply
apply these to each right hand side (or, alternatively, to form the LU factorizations,
and then perform triangular solves). The cost $T_{\rm inv}$ of computing the inverses,
and the cost $T_{\rm apply}$ of applying them then satisfy
\begin{align*}
T_{\rm inv}   &\sim N_{\rm P}^{3}\,N_{\rm G}^{3}\,N_{\rm F},\\
T_{\rm apply} &\sim N_{\rm P}^{2}\,N_{\rm G}^{2}\,N_{\rm F}.
\end{align*}

\lsp

We make some practical observations:
\begin{itemize}
\item The cost of forming the matrices by far dominates the other costs unless
the kernel is either smooth, or analytic formulas for $k_{n}$ are available.

\lsp

\item The scheme is highly efficient in situations where the same equation
needs to be solved for a sequence of different right hand sides. Given an additional right hand side,
the added cost $T_{\rm solve}$ is given by
$$
T_{\rm solve} \sim N_{\rm P}\,N_{\rm G}\,N_{\rm F}\,\log(N_{\rm F}) +
N_{\rm P}^{2}\,N_{\rm G}^{2}\,N_{\rm F},
$$
with a very small constant of proportionality. We note that this cost
remains small even if an analytic formula for $k_{n}$ is not
available.

\lsp

\item The system matrices $I + A_{n}$ often have internal structure that allow
them to be inverted using ``fast methods'' such as, \textit{e.g.}, those in
\cite{Martinsson:04a}. The cost of inversion and application can then be
accelerated to near optimal complexity.

\end{itemize}


\section{Simplifications for the double layer kernels associated with Laplace's equation}
\label{sec:laplace}


\subsection{The double layer kernels of Laplace's equation}
\label{sec:doublelayer}

Let $D \subseteq \mathbb{R}^{3}$ be a bounded domain whose boundary is
given by a smooth surface $\Gamma$, let $E = \bar{D}^{\rm c}$ denote the domain exterior to
$D$, and let $\nb$ and be the outward unit normal to $D$. Consider the interior and exterior Dirichlet problems
of potential theory \cite{Guenther:88a},
\begin{align}
    & \Delta u = 0 \hspace{0.5em} \textrm{in} \hspace{0.5em} D, \hspace{1em} u = f \hspace{0.5em} \textrm{on} \hspace{0.5em} \Gamma,
    \hspace{2em} \textrm{(interior Dirichlet problem)} \label{eq:intDir} \\
    & \Delta u = 0 \hspace{0.5em} \textrm{in} \hspace{0.5em} E, \hspace{1em} u = f \hspace{0.5em} \textrm{on} \hspace{0.5em} \Gamma.
    \hspace{2em} \textrm{(exterior Dirichlet problem)} \label{eq:extDir}
\end{align}
The solutions to (\ref{eq:intDir}) and (\ref{eq:extDir}) can be written in the respective forms
\begin{align*}
    & u(\xb) = \int_{\Gamma} \frac{\nb(\xb') \cdot (\xb - \xb')}{4 \pi |\xb - \xb'|^{3}} \sigma(\xb') \, dA(\xb'),
    \hspace{1em} \xb \in D,  \\
    & u(\xb) = \int_{\Gamma} \( -\frac{\nb(\xb') \cdot (\xb - \xb')}{4 \pi |\xb - \xb'|^{3}} + \frac{1}{4 \pi |\xb - \xb_{0}|} \) \sigma(\xb') \, dA(\xb'),
    \hspace{1em} \xb \in E, \hspace{0.5em} \xb_{0} \in D,
\end{align*}
where $\sigma$ is a boundary charge distribution that can be determined using the boundary
conditions. The resulting equations are
\begin{align}
    -\frac{1}{2} \sigma(\xb) &+ \int_{\Gamma} \frac{\nb(\xb') \cdot (\xb - \xb')}{4 \pi |\xb - \xb'|^{3}} \sigma(\xb') \, dA(\xb') = f(\xb),
    \hspace{1em} \xb \in \Gamma, \label{eq:intDirBie} \\
    -\frac{1}{2} \sigma(\xb) &+ \int_{\Gamma} \( -\frac{\nb(\xb') \cdot (\xb - \xb')}{4 \pi |\xb - \xb'|^{3}} + \frac{1}{4 \pi |\xb - \xb_{0}|} \) \sigma(\xb') \, dA(\xb') = f(\xb),
    \hspace{1em} \xb \in \Gamma. \label{eq:extDirBie}
\end{align}

\remark{There are other integral formulations for the solution to Laplace's equation.
The double layer formulation presented here is a good choice in that it provides an
integral operator that leads to well conditioned linear systems.  However, the
methodology of this paper is equally applicable to single-layer formulations that
lead to first kind Fredholm BIEs.}


\subsection{Separation of variables}

Using the procedure given in Section \ref{sec:sepVar}, if $\Gamma = \gamma \times \mathbb{T}$, then (\ref{eq:intDir}) and (\ref{eq:extDir}) can be recast as a series of BIEs defined along $\gamma$.  We express $\nb$ in cylindrical coordinates as
\begin{equation*}
    \nb(\xb') = (n_{r'} \cos \theta', n_{r'} \sin \theta', n_{z'}).
\end{equation*}
Further,
\begin{align*}
    |\xb - \xb'|^{2} &= (r \cos \theta - r' \cos \theta')^{2} + (r \sin \theta - r' \sin \theta')^{2} + (z - z')^2 \\
    &= r^{2} + (r')^{2} - 2 r r' (\sin \theta \sin \theta' + \cos \theta \cos \theta') + (z - z')^2 \\
    &= r^{2} + (r')^{2} - 2 r r' \cos(\theta - \theta') + (z - z')^{2}
\end{align*}
and
\begin{align*}
    \nb(\xb') \cdot (\xb - \xb') &= (n_{r'} \cos \theta', n_{r'} \sin \theta', n_{z'}) \cdot
    (r \cos \theta - r' \cos \theta', r \sin \theta - r' \sin \theta', z - z') \\
    &= n_{r'} r( \sin \theta \sin \theta' + \cos \theta \cos \theta') - n_{r'} r' + n_{z'}(z - z') \\
    &= n_{r'} (r \cos(\theta - \theta') - r') + n_{z'}(z - z').
\end{align*}
Then for a point $\xb' \in \Gamma$, the kernel of the internal Dirichlet problem can be expanded as
\begin{equation*}
    \frac{\nb(\xb') \cdot (\xb - \xb')}{4 \pi |\xb - \xb'|^{3}} = \frac{1}{\sqrt{2 \pi}} \sum_{n \in \mathbb{Z}} e^{i n (\theta - \theta')} d^{(i)}_{n}(r,z,r',z'),
\end{equation*}
where
\begin{equation*}
    d^{(i)}_{n}(r,z,r',z') = \frac{1}{\sqrt{32 \pi^{3}}} \int_{\mathbb{T}} e^{-i n \theta} \left[
     \frac{n_{r'} ( r \cos \theta - r' ) + n_{z'}(z - z')}
        {(r^{2} + (r')^{2} - 2 r r' \cos \theta + (z - z')^{2})^{3/2}} \right] \, d \theta.
\end{equation*}
Similarly, the kernel of the external Dirichlet problem can be written as
\begin{equation*}
    -\frac{\nb(\xb') \cdot (\xb - \xb')}{4 \pi |\xb - \xb'|^{3}} + \frac{1}{4 \pi | \xb - \xb_{0} |} =
    \frac{1}{\sqrt{2 \pi}} \sum_{n \in \mathbb{Z}} e^{i n (\theta - \theta')} d^{(e)}_{n}(r,z,r',z'),
\end{equation*}
with
\begin{align*}
    d^{(e)}_{n}(r,z,r',z') = \frac{1}{\sqrt{32 \pi^{3}}} & \int_{\mathbb{T}} e^{-i n \theta} \biggl(
     -\frac{n_{r'} ( r \cos \theta - r' ) + n_{z'}(z - z')}
        {(r^{2} + (r')^{2} - 2 r r' \cos \theta + (z - z')^{2})^{3/2}} + \\
    & + \frac{1}{(r^{2} + r_{0}^{2} - 2 r r_{0} \cos \theta + (z - z_{0})^{2})^{1/2}} \biggr)
     \, d \theta,
\end{align*}
where $\xb_{0}$ has been written in cylindrical coordinates as $(r_{0} \cos(\theta_{0}),r_{0} \sin(\theta_{0}),z_{0})$.
With the expansions of the kernels available, the procedure described in Section \ref{sec:genAlg} can be used to solve (\ref{eq:intDirBie}) and (\ref{eq:extDirBie}) by solving
\begin{equation}
    \sigma_{n}(r,z) + \sqrt{2 \pi} \int_{\gamma} d^{(i)}_{n}(r,r',z,z') \sigma_{n}(r',z') \, r' \, dl(r',z') = f_{n}(r,z)
    \label{eq:intDirFour}
\end{equation}
and
\begin{equation}
    \sigma_{n}(r,z) + \sqrt{2 \pi} \int_{\gamma} d^{(e)}_{n}(r,r',z,z') \sigma_{n}(r',z') \, r' \, dl(r',z') = f_{n}(r,z),
    \label{eq:extDirFour}
\end{equation}
respectively for $n = -N_{\rm{F}},-N_{\rm{F}}+1,\ldots,N_{\rm{F}}$.  Note that the kernels $d^{(i)}_{n}$ and $d^{(e)}_{n}$ contain a log-singularity when both $r' = r$ and $z' = z$.

Equivalently, (\ref{eq:intDirFour}) and (\ref{eq:extDirFour}) can be arrived at by considering Laplace's equation written in cylindrical coordinates,
\begin{equation*}
    \frac{\p^{2}u}{\p^{2}r} + \frac{1}{r}\frac{\p u}{\p r} +
    \frac{1}{r^{2}}\frac{\p^{2}u}{\p^{2} \theta} + \frac{\p^{2}u}{\p^{2}z}= 0,
\end{equation*}
Taking the Fourier transform of $u$ with respect to theta $\theta$ gives
\begin{equation*}
    \frac{\p^{2}u_{n}}{\p^{2}r} + \frac{1}{r}\frac{\p u_{n}}{\p r} -
    \frac{n^{2}}{r^{2}}\frac{\p^{2}u_{n}}{\p^{2} \theta} + \frac{\p^{2}u_{n}}{\p^{2}z}= 0, \hspace{1em} n \in \mathbb{Z},
\end{equation*}
where $e_{n} = e_{n}(\theta) = e^{i n \theta}/\sqrt{2\pi}$ and $u = \displaystyle\sum_{n \in \mathbb{Z}} e_{n} u_{n}$.  Then (\ref{eq:intDirFour}) and (\ref{eq:extDirFour}) are now associated with this sequence of PDEs.


\subsection{Evaluation of kernels}
\label{sec:recursion}

The values of $d^{(i)}_{n}$ and $d^{(e)}_{n}$ for $n = -N_{\rm{F}},-N_{\rm{F}}+1,\ldots,N_{\rm{F}}$ need to be computed efficiently and with high accuracy to construct the Nystr\"{o}m discretization of (\ref{eq:intDirFour}) and (\ref{eq:extDirFour}).  Note that the integrands of $d^{(i)}_{n}$ and $d^{(e)}_{n}$ are real valued and even functions on the interval $[-\pi,\pi]$.  Therefore, $d^{(i)}_{n}$ can be written as
\begin{equation}
    d_{n}^{(i)}(r,z,r',z') = \frac{1}{\sqrt{32 \pi^{3}}} \int_{\mathbb{T}} \left[ \frac{n_{r'} ( r \cos t - r' ) + n_{z'}(z - z')}
        {(r^{2} + (r')^{2} - 2 r r' \cos t + (z - z')^{2})^{3/2}} \right] \cos(n t) \, dt.
    \label{eq:doubleCoef}
\end{equation}
Note that $d^{(e)}_{n}$ can be written in a similar form.

This integrand is oscillatory and increasingly peaked at the origin as both $r' \rightarrow r$ and $z' \rightarrow z$.  As long as $r'$ and $r$ as well as $z'$ and $z$ are well separated, the integrand does not experience peaks near the origin, and as mentioned before, the FFT provides a fast and accurate way for calculating $d^{(i)}_{n}$ and $d^{(e)}_{n}$.

In regimes where the integrand is peaked, the FFT no longer provides a means of evaluating $d^{(i)}_{n}$ and $d^{(e)}_{n}$ with the desired accuracy.  One possible solution to this issue is applying adaptive quadrature to fully resolve the peak.  However, this must be done for each value of $n$ required and becomes prohibitively expensive if $N_{\rm{F}}$ is large.

Fortunately, an analytical solution to (\ref{eq:doubleCoef}) exists.  As noted in \cite{Cohl:99a}, the single-layer kernel can be expanded with respect to the azimuthal variable as
\begin{equation*}
    s(\xb,\xb') = \frac{1}{4 \pi | \xb - \xb' |} = \frac{1}{4 \pi^{2} \sqrt{r r'}} \sum_{n \in \mathbb{Z}} e^{i n (\theta - \theta')} \mathcal{Q}_{n-1/2}(\chi),
\end{equation*}
where $\mathcal{Q}_{n-1/2}$ is the half-integer degree Legendre function of the second kind and
\begin{equation*}
    \chi = \frac{r^{2} + (r')^{2} + (z - z')^{2}}{2 r r'}.
\end{equation*}
In light of this expansion, single-layer kernel can similarly be written as
\begin{align*}
    s(\xb,\xb') &= \frac{1}{4 \pi (r^{2} + (r')^{2} - 2 r r' \cos(\theta - \theta') + (z - z')^{2})^{1/2}} \\
    &= \frac{1}{\sqrt{2 \pi}} \sum_{n \in \mathbb{Z}} e^{i n (\theta - \theta')} s_{n}(r,z,r',z')
\end{align*}
where
\begin{align*}
    s_{n}(r,z,r',z') &= \frac{1}{\sqrt{32 \pi^{3}}} \int_{\mathbb{T}}
    \frac{\cos( n t)}{(r^{2} + (r')^{2} - 2 r r' \cos(t) + (z - z')^{2})^{1/2}} \, dt \\
    &= \frac{1}{ \sqrt{8 \pi^{3} r r'}} \int_{\mathbb{T}} \frac{\cos(nt)}{\sqrt{8 (\chi - \cos(t))} } \, dt \\
    &= \frac{1}{\sqrt{8 \pi^{3} r r'}} \mathcal{Q}_{n-1/2}(\chi).
\end{align*}

To find an analytical form for (\ref{eq:doubleCoef}), first note that in cylindrical coordinates the double-layer kernel can be written in terms of the single-layer kernel,
\begin{align*}
    \frac{\nb (\xb') \cdot (\xb - \xb')}{4 \pi |\xb - \xb'|^{3}} =& \frac{n_{r'} ( r \cos(\theta - \theta') - r' ) + n_{z'}(z - z')}
        {4 \pi (r^{2} + (r')^{2} - 2 r r' \cos(\theta - \theta') + (z - z')^{2})^{3/2}} \\
    =& \frac{1}{4 \pi} \biggl[ n_{r'} \frac{\p}{\p r'} \( \frac{1}{(r^{2} +
    (r')^{2} - 2 r r' \cos(\theta - \theta') + (z - z')^{2})^{1/2}} \) + \\
    &+ n_{z'} \frac{\p}{\p z'} \( \frac{1}{(r^{2} + (r')^{2} - 2 r r' \cos(\theta - \theta') + (z - z')^{2})^{1/2}} \) \biggr].
\end{align*}
The coefficients of the Fourier series expansion of the double-layer kernel are then given by $d^{(i)}_{n}$, which can be written using the previous equation as
\begin{align*}
    d_{n}^{(i)}(r,z,r',z') =& n_{r'} \int_{\mathbb{T}} \frac{\p}{\p r'} \( \frac{\cos(n t)}{ (32 \pi^{3} (r^{2} +
    (r')^{2} - 2 r r' \cos(t) + (z - z')^{2}))^{1/2}} \) \, dt + \\
    &+ n_{z'} \int_{\mathbb{T}} \frac{\p}{\p z'}
    \( \frac{\cos(n t)}{(32 \pi^{3}(r^{2} + (r')^{2} - 2 r r' \cos(t) + (z - z')^{2}))^{1/2}} \)  \, dt \\
    =& n_{r'} \frac{\p}{\p r'} \( \frac{1}{\sqrt{8 \pi^{3} r r'}} \mathcal{Q}_{n-1/2}(\chi) \) +
    n_{z'} \frac{\p}{\p z'} \( \frac{1}{\sqrt{8 \pi^{3} r r'}} \mathcal{Q}_{n-1/2}(\chi) \) \\
    =& \frac{1}{\sqrt{8 \pi^{3} r r'}}\[ n_{r'} \( \frac{\p \mathcal{Q}_{n-1/2}(\chi)}{\p \chi} \frac{\p \chi}{\p r'} -
    \frac{\mathcal{Q}_{n-1/2}(\chi)}{2 r'} \)
    + n_{z'} \frac{\p \mathcal{Q}_{n-1/2}(\chi)}{\p \chi} \frac{\p \chi}{\p z'} \].
\end{align*}
To utilize this form of $d^{(i)}_{n}$, set $\mu = \sqrt{\frac{2}{\chi + 1}}$ and note that
\begin{align*}
    &\frac{\p \chi}{\p r'} = \frac{(r')^{2} - r^{2} - (z - z')^{2}}{2 r (r')^{2}}, \\
    &\frac{\p \chi}{\p z'} = \frac{z' - z}{r r'}, \\
    &\mathcal{Q}_{-1/2}(\chi) = \mu K(\mu), \\
    &\mathcal{Q}_{1/2}(\chi) = \chi \mu K(\mu) - \sqrt{2(\chi + 1)} E(\mu), \\
    &\mathcal{Q}_{-n-1/2}(\chi) = \mathcal{Q}_{n-1/2}(\chi), \\
    &\mathcal{Q}_{n-1/2}(\chi) = 4 \frac{n-1}{2n-1} \chi \mathcal{Q}_{n-3/2}(\chi) - \frac{2n-3}{2n-1}\mathcal{Q}_{n-5/2}(\chi), \\
    &\frac{\p \mathcal{Q}_{n-1/2}(\chi)}{\p \chi} = \frac{2n-1}{2(\chi^2-1)} \( \chi \mathcal{Q}_{n-1/2} - \mathcal{Q}_{n-3/2} \),
\end{align*}
where $K$ and $E$ are the complete elliptic integrals of the first and second kinds, respectively.  The first two relations follow immediately from the definition of $\chi$ and the relations for the Legendre functions of the second kind can be found in \cite{Abramowitz:65a}.  With these relations in hand, the calculation of $d^{(i)}_{n}$ for $n = -N_{\rm{F}},-N_{\rm{F}}+1,\ldots,N_{\rm{F}}$ can be done accurately and efficiently when $r'$ and $r$ as well as $z'$ and $z$ are in close proximity.  The calculation of $d^{(e)}_{n}$ can be done analogously.

\remark{Note that the forward recursion relation for the Legendre functions $\mathcal{Q}_{n-1/2}(\chi)$ is unstable when $\chi > 1$.  In practice, the instability is mild when $\chi$ is near $1$ and can still be employed to accurately compute values in this regime.  Additionally, if stability becomes and issue, Miller's algorithm \cite{Gil:07a} can be used to calculate the values of the Legendre functions using the backwards recursion relation, which is stable for $\chi > 1$.}


\section{Numerical results} 
\label{sec:numRes}

\begin{figure}[!ht]
\begin{center}
\begin{minipage}{0.31\linewidth} \begin{flushleft}
\hspace{10mm}
\includegraphics[height=.75\linewidth]{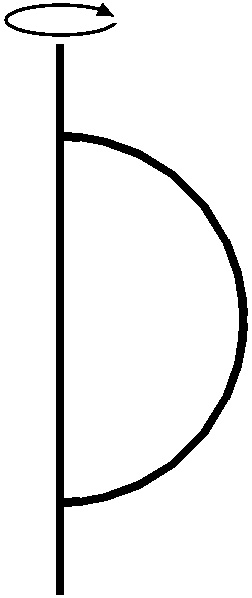} 
\end{flushleft} \end{minipage}
\begin{minipage}{0.31\linewidth} \begin{center}
\includegraphics[height =.7\linewidth]{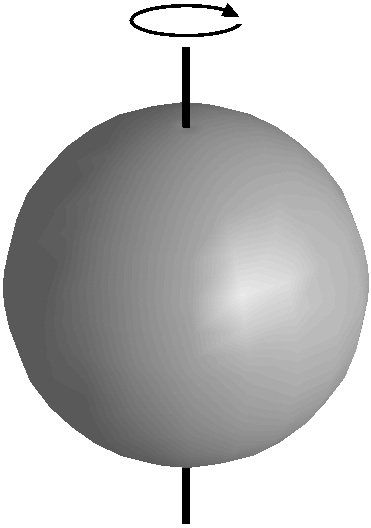}
\end{center} \end{minipage}
\\
\begin{minipage}{1\linewidth}\begin{center} (a) \end{center} \end{minipage}
\\ \vspace{4mm}
\begin{minipage}{0.31\linewidth} \begin{flushleft}
\hspace{10mm}
\includegraphics[height=.75\linewidth]{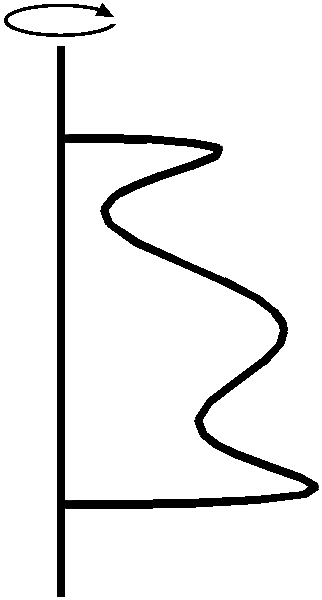} 
\end{flushleft} \end{minipage}
\begin{minipage}{0.31\linewidth} \begin{center}
\includegraphics[height =.65\linewidth]{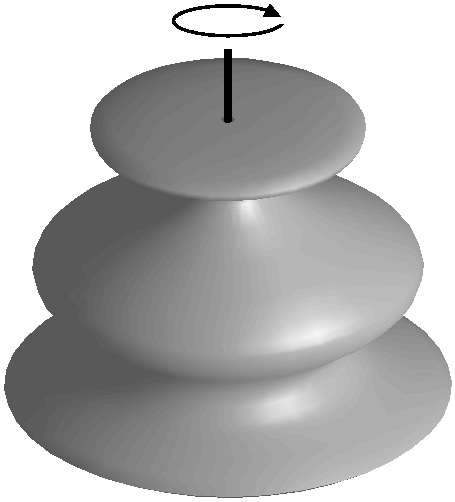}
\end{center} \end{minipage}
\\
\begin{minipage}{1\linewidth}\begin{center} (b) \end{center} \end{minipage}
\\ \vspace{4mm}
\begin{minipage}{0.31\linewidth} \begin{flushleft}
\hspace{10mm}
\includegraphics[height=.75\linewidth]{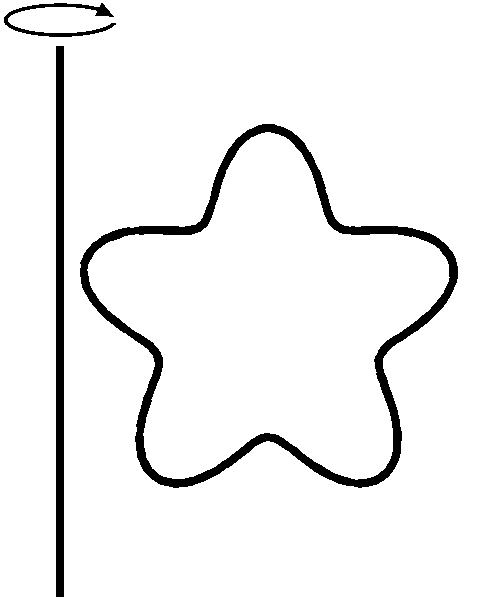} 
\end{flushleft} \end{minipage}
\begin{minipage}{0.31\linewidth} \begin{center}
\includegraphics[height =.55\linewidth]{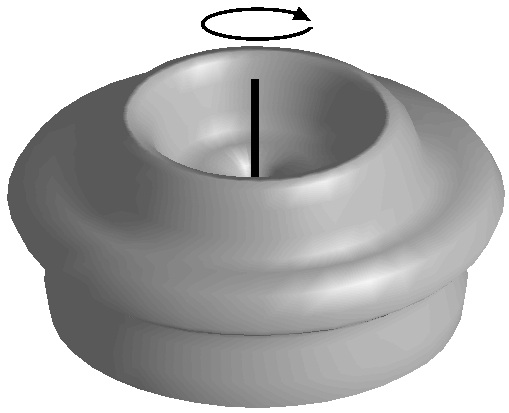}
\end{center} \end{minipage}
\\
\begin{minipage}{1\linewidth}\begin{center} (c) \end{center} \end{minipage}
\\
\caption{Domains used in numerical examples. All items are rotated about the vertical axis.  (a) A sphere.  (b) A wavy block.  (c) A starfish torus.  }
\label{fig:domains}
\end{center}
\end{figure}

This section describes several numerical experiments performed to assess
the efficiency and accuracy of the the numerical scheme outlined in
Section \ref{sec:summary}. All experiments were executed for the double
layer kernels associated with Laplace's equation, calculated using the
recursion relation described in Section \ref{sec:recursion}. Note that the kernels
in this case give us the property that $A_{-n} = A_{n}$, and so we need only to invert
$N_{\rm F} + 1$ matrices. Further, the FFT used here is complex-valued, and a real-valued FFT
would yield a significant decrease in computation time. The
geometries investigated are described in Figure \ref{fig:domains}. The
generating curves were parameterized by arc length, and split into
$N_{\rm P}$ equisized panels. A 10-point Gaussian quadrature has been
used along each panel, with the modified quadratures used to handle the
integrable singularities in the kernel.  These quadratures are listed in
the appendix.  The algorithm was implemented in MatLab and the
experiments were run on a 2.66GHz Intel Quad Core with 6Gb of RAM. All
timings were averaged over 10 runs.

\subsection{Computational costs}

Using the domain in Figure \ref{fig:domains}(a) and the interior Dirichlet problem, timing results are given in Table 1.
The reported results include:

\lsp

\begin{tabular}{l l}
$N_{\rm{P}}$ & the number  of panels used to discretize the contour \\
$N_{\rm{F}}$ & the Fourier truncation parameter (we keep $2N_{\rm F}+1$ modes) \\
$T_{\textrm{setup}}$ & time to setup the discrete system, excluding construction of the linear systems \\
$T_{\textrm{mat}}$ & time to construct the linear systems (utilizing the recursion relation) \\
$T_{\textrm{inv}}$ & time to invert the linear systems \\
$T_{\textrm{fft}}$ & time to Fourier transform the right hand side and the solution \\
$T_{\textrm{apply}}$ & time to apply the inverse to the right hand side \\
\end{tabular}

\lsp

\begin{figure}[!ht]
\begin{center}
\begin{minipage}{0.49\linewidth} \begin{center}
\includegraphics[width =.9\linewidth]{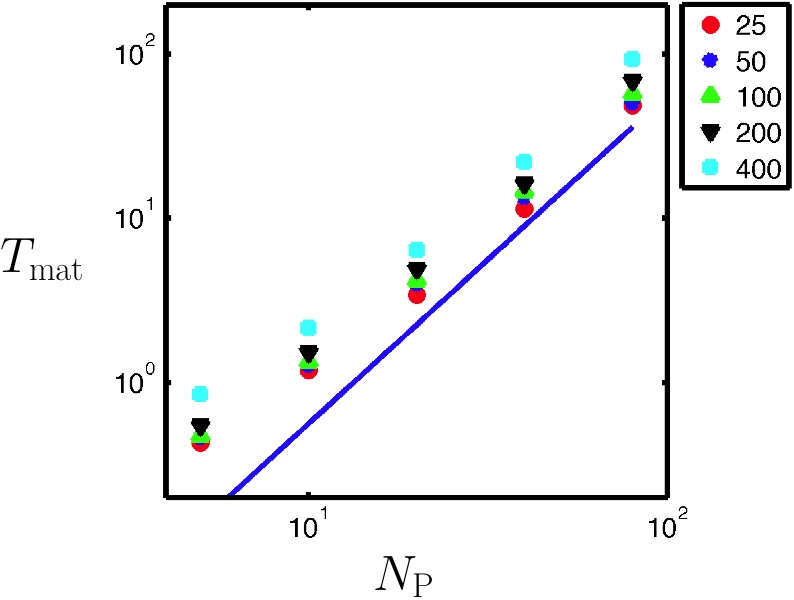}
\end{center} \end{minipage}
\begin{minipage}{0.49\linewidth} \begin{center}
\includegraphics[width =.9\linewidth]{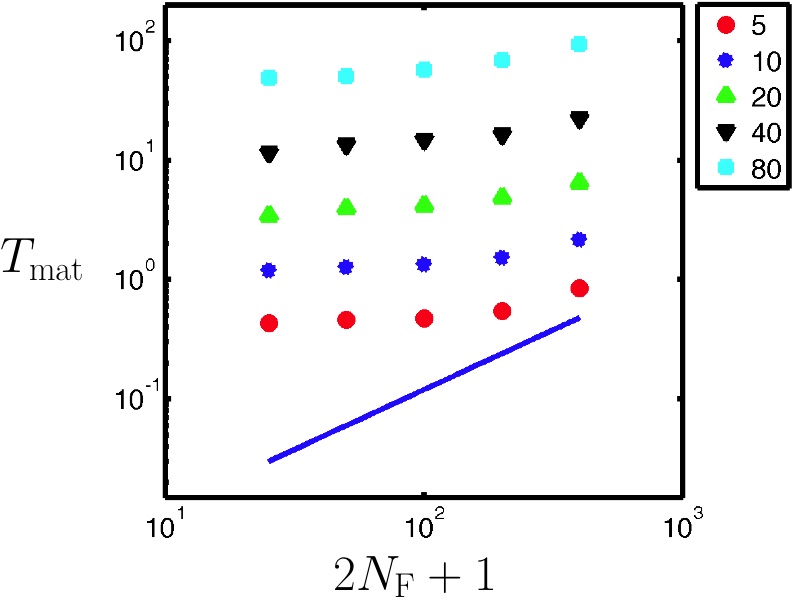}
\end{center} \end{minipage} \\ \vspace{2mm}
\begin{minipage}{0.48\linewidth}\begin{center} (a) \end{center} \end{minipage}
\begin{minipage}{0.48\linewidth}\begin{center} (b) \end{center} \end{minipage}
\vspace{-2mm}
\caption{Timings for the construction of the linear systems. (a) Scaling as the number of panels increases.  The solid line corresponds to growth with the number of panels squared.  The legend refers to the number of Fourier modes used.  (b) Scaling as the number of Fourier modes increases.  The solid line corresponds to linear growth.  The legend refers to the number of panels used.}
\label{fig:scalings}
\end{center}
\end{figure}

The most expensive component of the calculation is the construction of the linear systems.
This is primarily a result of the cost of evaluating the kernel and applying the modified
quadrature rules.  Table 2 compares the use of the recursion relation in evaluating the
kernel when it is near-singular to using an adaptive Gaussian quadrature.  The efficiency
of the recursion relation is clearly evident in this case.  Figure \ref{fig:scalings} plots
the time to construct the linear systems as the number of panels and as the number of Fourier
modes increases.  The time to construct the systems as the number of panels increases grows
as $O(N_{\rm{P}}^{2} N_{\rm{G}}^{2})$ as expected, and the timings go as $O(N_{\rm{F}})$
as the number of Fourier modes increase.  Note that it is rather inexpensive to increase the 
number of Fourier modes used, making the scheme
particularly attractive in environments where a large number is required.  
As an artifact of the small lengths of the vectors being Fourier
transformed, the asymptotic behavior of the timings as the number of Fourier modes increases 
is not evident until we reach a large number of modes.

We observe that the largest problem reported in Table 1 involves
$320\,000$ degrees of freedom. The method requires $2.2$ minutes of
pre-computation for this example, and is then capable of computing a
solution $u$ from a given data function $f$ in $0.46$ seconds.

\begin{table}[!ht]
\begin{center}
\begin{tabular}{| c | c | c | c | c | c | c | }
\hline
$N_{\rm{P}}$ & $2N_{\rm{F}} + 1$ & $T_{\textrm{setup}}$ & $T_{\textrm{mat}}$ &
$T_{\textrm{inv}}$ & $T_\textrm{{fft}}$ & $T_{\textrm{apply}}$ \\
\hline
5 & 25 & 9.10e-02 & 4.31e-01 & 2.48e-03 & 1.17e-03 & 3.82e-04 \\
10 & 25 & 1.10e-01 & 1.19e+00 & 9.31e-03 & 1.93e-03 & 5.21e-04 \\
20 & 25 & 1.22e-01 & 3.42e+00 & 4.47e-02 & 3.42e-03 & 1.20e-03 \\
40 & 25 & 1.68e-01 & 1.15e+01 & 2.90e-01 & 6.66e-03 & 4.94e-03 \\
80 & 25 & 2.60e-01 & 4.90e+01 & 2.12e+00 & 1.32e-02 & 1.94e-02 \\ \hline
5 & 50 & 9.63e-02 & 4.59e-01 & 4.79e-03 & 1.63e-03 & 6.85e-04 \\
10 & 50 & 1.16e-01 & 1.27e+00 & 1.72e-02 & 2.81e-03 & 9.71e-04 \\
20 & 50 & 1.54e-01 & 3.99e+00 & 8.89e-02 & 5.47e-03 & 2.86e-03 \\
40 & 50 & 2.32e-01 & 1.34e+01 & 5.80e-01 & 1.08e-02 & 1.00e-02 \\
80 & 50 & 4.04e-01 & 5.05e+01 & 4.27e+00 & 2.13e-02 & 3.89e-02 \\ \hline
5 & 100 & 1.13e-01 & 4.72e-01 & 9.20e-03 & 2.74e-03 & 1.29e-03 \\
10 & 100 & 1.49e-01 & 1.34e+00 & 3.34e-02 & 4.98e-03 & 1.99e-03 \\
20 & 100 & 2.20e-01 & 4.19e+00 & 1.75e-01 & 9.80e-03 & 6.20e-03 \\
40 & 100 & 3.75e-01 & 1.47e+01 & 1.16e+00 & 2.06e-02 & 2.01e-02 \\
80 & 100 & 6.53e-01 & 5.76e+01 & 8.34e+00 & 4.13e-02 & 7.99e-02 \\ \hline
5 & 200 & 1.48e-01 & 5.46e-01 & 1.92e-02 & 5.14e-03 & 2.68e-03 \\
10 & 200 & 2.24e-01 & 1.53e+00 & 6.89e-02 & 9.53e-03 & 4.48e-03 \\
20 & 200 & 3.59e-01 & 4.90e+00 & 3.44e-01 & 1.88e-02 & 1.25e-02 \\
40 & 200 & 6.45e-01 & 1.63e+01 & 2.29e+00 & 3.74e-02 & 4.13e-02 \\
80 & 200 & 1.22e+00 & 6.85e+01 & 1.68e+01 & 7.67e-02 & 1.58e-01 \\ \hline
5 & 400 & 2.11e-01 & 8.49e-01 & 3.81e-02 & 1.01e-02 & 5.46e-03 \\
10 & 400 & 3.45e-01 & 2.15e+00 & 1.29e-01 & 1.92e-02 & 9.59e-03 \\
20 & 400 & 5.83e-01 & 6.47e+00 & 6.85e-01 & 3.81e-02 & 2.54e-02 \\
40 & 400 & 1.09e+00 & 2.22e+01 & 4.49e+00 & 7.60e-02 & 7.94e-02 \\
80 & 400 & 2.27e+00 & 9.33e+01 & 3.27e+01 & 1.51e-01 & 3.06e-01 \\ \hline
\end{tabular}
\vspace{2mm}
\caption{Timing results in seconds performed for the domain given in Figure \ref{fig:domains}(a) for the interior Dirichlet problem.}
\end{center}
\end{table}

\begin{table}[!ht]
\begin{center}
\begin{tabular}{| c | c | c |}
\hline
$2N_{\rm{F}} + 1$ & Composite Quadrature & Recursion Relation \\
\hline
25 & 1.9 & 0.43 \\
50 & 3.1 & 0.46 \\
100 & 6.6 & 0.47 \\
200 & 18.9 & 0.55 \\
\hline
\end{tabular}
\vspace{2mm}
\caption{Timing comparison in seconds for constructing the matrices $(I + A_{n})$ using composite Gaussian quadrature and the recursion relation described in Section \ref{sec:recursion}  to evaluate $k_{n}$ for diagonal and near diagonal blocks.  The FFT is used to evaluate $k_{n}$ at all other entries.  $2N_{\rm{F}} + 1$ is the total number of Fourier modes used. 5 panels were used to discretize the boundary.}
\end{center}
\end{table}

\subsection{Accuracy and conditioning of discretization}

The accuracy of the discretization has been tested using the interior and exterior Dirichlet
problems on the domains given in Figure \ref{fig:domains}.  Exact solutions were generated by
placing a few random point charges outside of the domain where the solution was calculated.
The solution was evaluated at points defined on a sphere encompassing (or interior to) the
boundary.  The errors reported in Tables 3, 4, and 5 are relative errors measured in the
$l^{\infty}$-norm, $||u_{\epsilon} - u||_{\infty} / || u ||_{\infty}$, where $u$ is the
exact potential and $u_{\epsilon}$ is the potential obtained from the numerical solution.

In all cases, 10 digits of accuracy has been obtained from a discretization involving a
relatively small number of panels, due to the rapid convergence of the Gaussian quadrature.
This is especially advantageous, as the most expensive component of the algorithm is the
construction of the linear systems, the majority of the cost being directly related to
the number of panels used.  Further, the number of Fourier modes required to obtain 10
digits of accuracy is on the order of 100 modes.  Although not investigated here, the
discretization technique naturally lends itself to nonuniform refinement of the surface,
allowing one to resolve features of the surface that require finer resolution.

The number of correct digits obtained as the number of panels and number of Fourier modes
increases eventually stalls.  This is a result of a loss of precision in determining the
kernels, as well as cancelation errors incurred when evaluating interactions between nearby
points.  This is especially prominent with the use of Gaussian quadratures, as points cluster
near the ends of the panels.  If more digits are required, high precision arithmetic can be
employed in the setup phase of the algorithm.

\begin{table}[!ht]
\begin{center}
\begin{minipage}{0.99\linewidth} \begin{center}
 \begin{tabular}{| c | c | c | c | c | c |}
 \hline
 \multicolumn{1}{|c|}{$N_{\rm P}$} & \multicolumn{5}{c|}{$2N_{\rm F} + 1$} \\
 \hline
- &25&50&100&200&400    \\
   \hline
5&1.93869e-04&4.10935e-07&5.37883e-08&5.37880e-08&5.37880e-08 \\
10&1.93869e-04&4.10513e-07&3.27169e-12&6.72270e-13&6.72270e-13 \\
20&1.93869e-04&4.10513e-07&3.30601e-12&1.66132e-13&1.66132e-13 \\
40&1.93869e-04&4.10513e-07&3.23162e-12&8.28568e-14&8.28568e-14 \\
80&1.93869e-04&4.10512e-07&2.92918e-12&2.92091e-13&2.92091e-13 \\
  \hline
\end{tabular}
\end{center} \end{minipage} \\
\vspace{2mm}
\caption{Error in internal Dirichlet problem solved on domain (a) in Figure \ref{fig:domains}.}
\end{center}
\end{table}

\begin{table}[!ht]
\begin{center}
\begin{minipage}{0.99\linewidth} \begin{center}
 \begin{tabular}{| c | c | c | c | c | c |}
 \hline
 \multicolumn{1}{|c|}{$N_{\rm P}$} & \multicolumn{5}{c|}{$2N_{\rm F} + 1$} \\
 \hline
- &25&50&100&200&400    \\
   \hline
5&9.11452e-04&9.11464e-04&9.11464e-04&9.11464e-04&9.11464e-04 \\
10&4.15377e-05&4.15416e-05&4.15416e-05&4.15416e-05&4.15416e-05 \\
20&6.31923e-07&1.29234e-07&1.29235e-07&1.29235e-07&1.29235e-07 \\
40&7.04741e-07&3.10049e-11&3.08152e-11&3.08305e-11&3.08359e-11 \\
80&7.04779e-07&5.62558e-11&5.05306e-11&5.05257e-11&5.05232e-11 \\
  \hline
\end{tabular}
\end{center} \end{minipage} \\
\vspace{2mm}
\caption{Error in external Dirichlet problem solved on domain (b) in Figure \ref{fig:domains}.}
\end{center}
\end{table}

\begin{table}[!ht]
\begin{center}
\begin{minipage}{0.99\linewidth} \begin{center}
 \begin{tabular}{| c | c | c | c | c | c |}
 \hline
 \multicolumn{1}{|c|}{$N_{\rm P}$} & \multicolumn{5}{c|}{$2N_{\rm F} + 1$} \\
 \hline
- &25&50&100&200&400    \\
   \hline
5&3.80837e-04 & 3.83707e-04 & 3.83707e-04 & 3.83707e-04 & 3.83707e-04 \\
10&2.41602e-05 & 6.81564e-06 & 6.81556e-06 & 6.81556e-06 & 6.81556e-06 \\
20&3.03272e-05 & 5.98506e-09 & 2.53980e-11 & 2.54112e-11 & 2.54118e-11 \\
40&3.03272e-05 & 6.01273e-09 & 6.95662e-12 & 6.94592e-12 & 6.94546e-12 \\
80&3.03272e-05 & 6.01059e-09 & 5.25217e-12 & 5.26674e-12 & 5.26515e-12 \\
  \hline
\end{tabular}
\end{center} \end{minipage} \\
\vspace{2mm}
\caption{Error in external Dirichlet problem solved on domain (c) in Figure \ref{fig:domains}.}
\label{tbl:InterpAcc}
\end{center}
\end{table}

\begin{figure}[htbp] 
   \centering
   \includegraphics[width=0.5\linewidth]{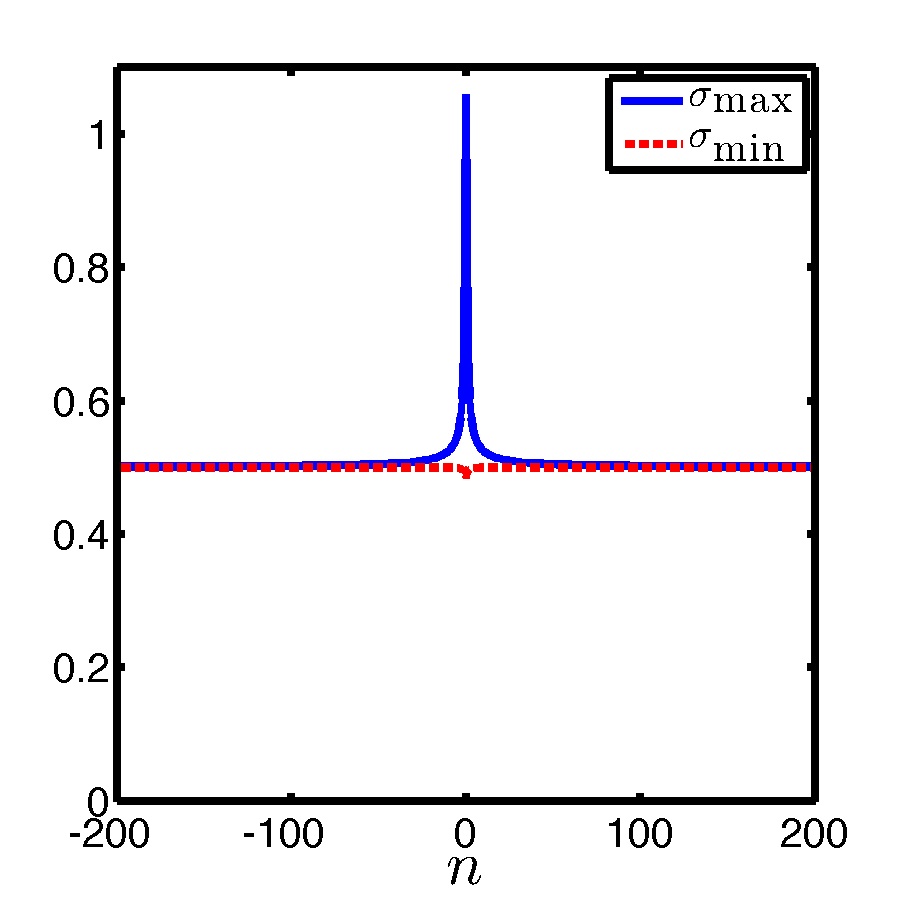} 
\caption{Maximum and minimum singular values for the matrices resulting from an 80 panel discretization of a sphere using 400 Fourier modes, where $n$ is the the matrix associated with the $n^{\textrm{th}}$ Fourier mode.}
\label{fig:cond}
\end{figure}

Figure \ref{fig:cond} shows the singular values as well as the condition numbers of an 80 panel discretization for the $N_{\rm{F}} = -200,\ldots,200$ Fourier modes used in the discretization of the interior Dirichlet problem, on the domain shown in Figure \ref{fig:domains}(a).  The integral equations of this paper are second kind Fredholm equations, and generally lead to well-conditioned systems.  As seen in Figure \ref{fig:cond}, this hold true for the discretization presented in this paper.

\section{Generalizations and Conclusions}
\label{sec:conclusions}

This paper describes a numerical technique for computing solutions to
boundary integral equations defined on axisymmetric surfaces in
$\mathbb{R}^{3}$ with no assumption on the loads being axisymmetric. The
technique is introduced as a generic method with only very mild
conditions imposed on the kernel; specifically, we assume that the kernel
has an integrable singularity at the diagonal, and that it is rotationally symmetric (in the sense that
(\ref{eq:kernel_condition}) holds). The technique described improves upon
previous work in two regards:
\begin{enumerate}
\item A highly accurate quadrature scheme for kernels with integrable
      singularities is introduced. Numerical experiments indicate that solutions with
      a relative accuracy of $10^{-10}$ or better can easily be constructed.
\item A rapid technique for numerically constructing the kernel functions
      $k_{n}$ in (\ref{eq:fred2}) is introduced. It works when $k$ is either
      the single or the double layer potential associated with Laplace's equation.
      The technique is a hybrid scheme that relies on the FFT when possible, and
      uses recursion relations for Legendre functions when not. The resulting
      scheme is fast enough that a problem involving $320\,000$ degrees of freedom
      can be solved in 2.2 minutes on a standard desktop PC. Once one problem has
      been solved, additional right hand sides can be processed in 0.46 seconds.
\end{enumerate}

Some possible extensions of this work include: (1) Acceleration of the
solution of the linear systems using fast methods, (2) extension of
recursion relation to the Helmholtz equation, and (3) extension of the
algorithm to problems involving multiple bodies whose axes of symmetry
are not necessarily aligned.

\lsp

\noindent
\textbf{Acknowledgements:} The authors are grateful for the support of NSF Grant DMS-0602284 (P. M. Young) and NSF Grant DMS-0748488 (P. G. Martinsson).  The authors have greatly benefited from valuable
suggestions made by Vladimir Rokhlin of Yale University.

\bibliographystyle{plain}
\bibliography{MyReferences}

\clearpage

\section*{APPENDIX OF QUADRATURE NODES AND WEIGHTS}

\footnotesize

\begin{center}
\begin{minipage}{0.48\linewidth} \begin{center}
\begin{tabular}{ | c | c | }
\hline
  \multicolumn{2}{|c|}{10 Point Gauss-Legendre Rule for}   \\
  \multicolumn{2}{|c|}{integrals of the form $\int_{-1}^{1}f(x) \, dx$}   \\
\hline
    NODES & WEIGHTS  \\
\hline
    -9.739065285171716e-01  & 6.667134430868814e-02 \\
    -8.650633666889845e-01  & 1.494513491505806e-01 \\
    -6.794095682990244e-01  & 2.190863625159820e-01 \\
    -4.333953941292472e-01  & 2.692667193099963e-01 \\
    -1.488743389816312e-01  & 2.955242247147529e-01 \\
     1.488743389816312e-01  & 2.955242247147529e-01 \\
     4.333953941292472e-01  & 2.692667193099963e-01 \\
     6.794095682990244e-01  & 2.190863625159820e-01 \\
     8.650633666889845e-01  & 1.494513491505806e-01 \\
     9.739065285171716e-01  & 6.667134430868814e-02 \\
\hline
\end{tabular}
\end{center} \end{minipage}
\begin{minipage}{0.48\linewidth} \begin{center}
\begin{tabular}{ | c | c | }
\hline
  \multicolumn{2}{|c|}{20 point quadrature rule for integrals}   \\
  \multicolumn{2}{|c|}{of the form $\int_{-1}^{1}f(x) + g(x) \log|x_{1} - x| \, dx$,} \\
  \multicolumn{2}{|c|}{where $x_{1}$ is a Gauss-Legendre node}   \\
\hline
  NODES & WEIGHTS  \\
\hline
    -9.981629455677877e-01     & 4.550772157144354e-03 \\
    -9.915520723139890e-01     & 8.062764683328619e-03 \\
    -9.832812993252168e-01     & 7.845621096866406e-03 \\
    -9.767801773920733e-01     & 4.375212351185101e-03 \\
    -9.717169387169078e-01     & 1.021414662954223e-02 \\
    -9.510630103726074e-01     & 3.157199356768625e-02 \\
    -9.075765988474132e-01     & 5.592493151946541e-02 \\
    -8.382582352569804e-01     & 8.310260847601852e-02 \\
    -7.408522006801963e-01     & 1.118164522164500e-01 \\
    -6.147619568252419e-01     & 1.401105427713687e-01 \\
    -4.615244999958006e-01     & 1.657233639623953e-01 \\
    -2.849772954295424e-01     & 1.863566566231937e-01 \\
    -9.117593460489747e-02     & 1.999093145144455e-01 \\
     1.119089520342051e-01     & 2.046841584582030e-01 \\
     3.148842536644393e-01     & 1.995580161940930e-01 \\
     5.075733846631832e-01     & 1.841025430283230e-01 \\
     6.797470718157004e-01     & 1.586456191174843e-01 \\
     8.218833662202629e-01     & 1.242680229936124e-01 \\
     9.258924858821892e-01     & 8.273794370795576e-02 \\
     9.857595961761246e-01     & 3.643931593123844e-02 \\
\hline
\end{tabular}
\end{center} \end{minipage}\\ \vspace{5em}
\begin{minipage}{0.48\linewidth} \begin{center}
\begin{tabular}{ | c | c | }
\hline
  \multicolumn{2}{|c|}{20 point quadrature rule for integrals}   \\
  \multicolumn{2}{|c|}{of the form $\int_{-1}^{1}f(x) + g(x) \log|x_{2} - x| \, dx$,} \\
  \multicolumn{2}{|c|}{where $x_{2}$ is a Gauss-Legendre node}   \\
\hline
  NODES & WEIGHTS  \\
\hline
    -9.954896691005256e-01     & 1.141744473788874e-02 \\
    -9.775532683688947e-01     & 2.368593568061651e-02 \\
    -9.500346715183706e-01     & 3.027205199814611e-02 \\
    -9.192373372373420e-01     & 3.021809354380292e-02 \\
    -8.916563772395616e-01     & 2.397183723558556e-02 \\
    -8.727728136507039e-01     & 1.253574079839078e-02 \\
    -8.607963163061316e-01     & 2.070840476545303e-02 \\
    -8.201318720954396e-01     & 6.080709508468810e-02 \\
    -7.394732321355052e-01     & 1.002402801599464e-01 \\
    -6.204853512352519e-01     & 1.371499151597280e-01 \\
    -4.667290485167077e-01     & 1.693838059093582e-01 \\
    -2.840823320902124e-01     & 1.945292086962893e-01 \\
    -8.079364608026202e-02     & 2.103223087093422e-01 \\
     1.328455136645940e-01     & 2.149900928447852e-01 \\
     3.451233500669768e-01     & 2.074984762344433e-01 \\
     5.437321547508867e-01     & 1.877085225595498e-01 \\
     7.167077216635750e-01     & 1.564543949958065e-01 \\
     8.534299232009863e-01     & 1.156104890379952e-01 \\
     9.458275339169444e-01     & 6.859369195724087e-02 \\
     9.912353127269481e-01     & 2.390220989094312e-02 \\
\hline
\end{tabular}
\end{center} \end{minipage}
\begin{minipage}{0.48\linewidth} \begin{center}
\begin{tabular}{ | c | c | }
\hline
  \multicolumn{2}{|c|}{20 point quadrature rule for integrals}   \\
  \multicolumn{2}{|c|}{of the form $\int_{-1}^{1}f(x) + g(x) \log|x_{3} - x| \, dx$,} \\
  \multicolumn{2}{|c|}{where $x_{3}$ is a Gauss-Legendre node}   \\
\hline
  NODES & WEIGHTS  \\
\hline
    -9.930122613589740e-01     & 1.779185041193254e-02 \\
    -9.643941806993207e-01     & 3.870503119897836e-02 \\
    -9.175869559770760e-01     & 5.371120494602663e-02 \\
    -8.596474181980754e-01     & 6.073467932536858e-02 \\
    -7.990442708271941e-01     & 5.901993373645797e-02 \\
    -7.443700671611690e-01     & 4.905519963921684e-02 \\
    -7.031684479828371e-01     & 3.249237036645046e-02 \\
    -6.811221147275545e-01     & 1.335394660596527e-02 \\
    -6.579449960254029e-01     & 4.151626407911676e-02 \\
    -5.949471688137100e-01     & 8.451456165895121e-02 \\
    -4.893032793226841e-01     & 1.262522607368499e-01 \\
    -3.441659232382107e-01     & 1.628408264966550e-01 \\
    -1.665388322404095e-01     & 1.907085686614375e-01 \\
     3.344207582228461e-02     & 2.071802230953481e-01 \\
     2.434356263087524e-01     & 2.105274833603497e-01 \\
     4.498696863725133e-01     & 2.000282912446872e-01 \\
     6.389777518528792e-01     & 1.760212445284564e-01 \\
     7.978632877793501e-01     & 1.399000904426490e-01 \\
     9.155180703268415e-01     & 9.402669072995991e-02 \\
     9.837258757826489e-01     & 4.161927873514264e-02 \\
\hline
\end{tabular}
\end{center} \end{minipage}\\
\end{center}

\begin{center}
\begin{minipage}{0.48\linewidth} \begin{center}
\begin{tabular}{ | c | c | }
\hline
  \multicolumn{2}{|c|}{20 point quadrature rule for integrals}   \\
  \multicolumn{2}{|c|}{of the form $\int_{-1}^{1}f(x) + g(x) \log|x_{4} - x| \, dx$,} \\
  \multicolumn{2}{|c|}{where $x_{4}$ is a Gauss-Legendre node}   \\
\hline
  NODES & WEIGHTS  \\
\hline
    -9.903478871133073e-01     & 2.462513260640712e-02 \\
    -9.504025146897784e-01     & 5.449201732062665e-02 \\
    -8.834986023815121e-01     & 7.799498604905293e-02 \\
    -7.974523551287549e-01     & 9.241688894090601e-02 \\
    -7.022255002503461e-01     & 9.619882322938848e-02 \\
    -6.087194789244920e-01     & 8.902783806614303e-02 \\
    -5.275278952351541e-01     & 7.181973054766198e-02 \\
    -4.677586540799037e-01     & 4.663017060126023e-02 \\
    -4.360689210457623e-01     & 1.794303974050253e-02 \\
    -4.121945474875853e-01     & 4.061799823415495e-02 \\
    -3.494226766911471e-01     & 8.507517518447759e-02 \\
    -2.425993523586304e-01     & 1.277525783357134e-01 \\
    -9.646839923908594e-02     & 1.628510773009247e-01 \\
     7.921243716767302e-02     & 1.863323765408308e-01 \\
     2.715178194484646e-01     & 1.958227701927855e-01 \\
     4.658440358656903e-01     & 1.903138548150517e-01 \\
     6.472213975763533e-01     & 1.700731513381802e-01 \\
     8.015601619414859e-01     & 1.365784674773513e-01 \\
     9.168056007307982e-01     & 9.239595239693155e-02 \\
     9.839468743284722e-01     & 4.103797108164931e-02 \\
\hline
\end{tabular}
\end{center} \end{minipage}
\begin{minipage}{0.48\linewidth} \begin{center}
\begin{tabular}{ | c | c | }
\hline
  \multicolumn{2}{|c|}{20 point quadrature rule for integrals}   \\
  \multicolumn{2}{|c|}{of the form $\int_{-1}^{1}f(x) + g(x) \log|x_{5} - x| \, dx$,} \\
  \multicolumn{2}{|c|}{where $x_{5}$ is a Gauss-Legendre node}   \\
\hline
  NODES & WEIGHTS  \\
\hline
    -9.883561797860961e-01     & 2.974603958509255e-02 \\
    -9.398305159297058e-01     & 6.657945456889164e-02 \\
    -8.572399919019390e-01     & 9.731775484182564e-02 \\
    -7.482086250804679e-01     & 1.190433988432928e-01 \\
    -6.228514167093102e-01     & 1.297088242013777e-01 \\
    -4.928317114329241e-01     & 1.282900896966494e-01 \\
    -3.702771193724617e-01     & 1.148917968875341e-01 \\
    -2.666412108172461e-01     & 9.074932908233864e-02 \\
    -1.916083010783277e-01     & 5.818196361216740e-02 \\
    -1.521937160593461e-01     & 2.224697059733435e-02 \\
    -1.233125650067164e-01     & 4.788826761346366e-02 \\
    -5.257959675044444e-02     & 9.237500180593534e-02 \\
     5.877314311857769e-02     & 1.287410543031414e-01 \\
     2.012559739993003e-01     & 1.541960911507042e-01 \\
     3.627988191760868e-01     & 1.665885274544506e-01 \\
     5.297121321076323e-01     & 1.648585116745725e-01 \\
     6.878399330187783e-01     & 1.491408089644010e-01 \\
     8.237603202215137e-01     & 1.207592726093190e-01 \\
     9.259297297557394e-01     & 8.212177982524418e-02 \\
     9.856881498392895e-01     & 3.657506268226379e-02 \\
\hline
\end{tabular}
\end{center} \end{minipage}\\ \vspace{5em}
\begin{minipage}{0.48\linewidth} \begin{center}
\begin{tabular}{ | c | c | }
\hline
  \multicolumn{2}{|c|}{20 point quadrature rule for integrals}   \\
  \multicolumn{2}{|c|}{of the form $\int_{-1}^{1}f(x) + g(x) \log|x_6 - x| \, dx$,} \\
  \multicolumn{2}{|c|}{where $x_{6}$ is a Gauss-Legendre node}   \\
\hline
  NODES & WEIGHTS  \\
\hline
    -9.856881498392895e-01     & 3.657506268226379e-02 \\
    -9.259297297557394e-01     & 8.212177982524418e-02 \\
    -8.237603202215137e-01     & 1.207592726093190e-01 \\
    -6.878399330187783e-01     & 1.491408089644010e-01 \\
    -5.297121321076323e-01     & 1.648585116745725e-01 \\
    -3.627988191760868e-01     & 1.665885274544506e-01 \\
    -2.012559739993003e-01     & 1.541960911507042e-01 \\
    -5.877314311857769e-02     & 1.287410543031414e-01 \\
     5.257959675044444e-02     & 9.237500180593534e-02 \\
     1.233125650067164e-01     & 4.788826761346366e-02 \\
     1.521937160593461e-01     & 2.224697059733435e-02 \\
     1.916083010783277e-01     & 5.818196361216740e-02 \\
     2.666412108172461e-01     & 9.074932908233864e-02 \\
     3.702771193724617e-01     & 1.148917968875341e-01 \\
     4.928317114329241e-01     & 1.282900896966494e-01 \\
     6.228514167093102e-01     & 1.297088242013777e-01 \\
     7.482086250804679e-01     & 1.190433988432928e-01 \\
     8.572399919019390e-01     & 9.731775484182564e-02 \\
     9.398305159297058e-01     & 6.657945456889164e-02 \\
     9.883561797860961e-01     & 2.974603958509255e-02 \\
\hline
\end{tabular}
\end{center} \end{minipage}
\begin{minipage}{0.48\linewidth} \begin{center}
\begin{tabular}{ | c | c | }
\hline
  \multicolumn{2}{|c|}{20 point quadrature rule for integrals}   \\
  \multicolumn{2}{|c|}{of the form $\int_{-1}^{1}f(x) + g(x) \log|x_{7} - x| \, dx$,} \\
  \multicolumn{2}{|c|}{where $x_{7}$ is a Gauss-Legendre node}   \\
\hline
  NODES & WEIGHTS  \\
\hline
    -9.839468743284722e-01     & 4.103797108164931e-02  \\
    -9.168056007307982e-01 &     9.239595239693155e-02  \\
    -8.015601619414859e-01 &     1.365784674773513e-01 \\
    -6.472213975763533e-01     & 1.700731513381802e-01 \\
    -4.658440358656903e-01     & 1.903138548150517e-01 \\
    -2.715178194484646e-01     & 1.958227701927855e-01 \\
    -7.921243716767302e-02     & 1.863323765408308e-01 \\
     9.646839923908594e-02     & 1.628510773009247e-01 \\
     2.425993523586304e-01     & 1.277525783357134e-01 \\
     3.494226766911471e-01     & 8.507517518447759e-02 \\
     4.121945474875853e-01     & 4.061799823415495e-02 \\
     4.360689210457623e-01     & 1.794303974050253e-02 \\
     4.677586540799037e-01     & 4.663017060126023e-02 \\
     5.275278952351541e-01     & 7.181973054766198e-02 \\
     6.087194789244920e-01     & 8.902783806614303e-02 \\
     7.022255002503461e-01     & 9.619882322938848e-02 \\
     7.974523551287549e-01     & 9.241688894090601e-02 \\
     8.834986023815121e-01     & 7.799498604905293e-02 \\
     9.504025146897784e-01     & 5.449201732062665e-02 \\
     9.903478871133073e-01     & 2.462513260640712e-02 \\
\hline
\end{tabular}
\end{center} \end{minipage}\\
\end{center}

\begin{center}
\begin{minipage}{0.48\linewidth} \begin{center}
\begin{tabular}{ | c | c | }
\hline
  \multicolumn{2}{|c|}{20 point quadrature rule for integrals}   \\
  \multicolumn{2}{|c|}{of the form $\int_{-1}^{1}f(x) + g(x) \log|x_{8} - x| \, dx$,} \\
  \multicolumn{2}{|c|}{where $x_{8}$ is a Gauss-Legendre node}   \\
\hline
  NODES & WEIGHTS  \\
\hline
    -9.837258757826489e-01     & 4.161927873514264e-02 \\
    -9.155180703268415e-01     & 9.402669072995991e-02 \\
    -7.978632877793501e-01     & 1.399000904426490e-01 \\
    -6.389777518528792e-01     & 1.760212445284564e-01 \\
    -4.498696863725133e-01     & 2.000282912446872e-01 \\
    -2.434356263087524e-01     & 2.105274833603497e-01 \\
    -3.344207582228461e-02     & 2.071802230953481e-01 \\
     1.665388322404095e-01     & 1.907085686614375e-01 \\
     3.441659232382107e-01     & 1.628408264966550e-01 \\
     4.893032793226841e-01     & 1.262522607368499e-01 \\
     5.949471688137100e-01     & 8.451456165895121e-02 \\
     6.579449960254029e-01     & 4.151626407911676e-02 \\
     6.811221147275545e-01     & 1.335394660596527e-02 \\
     7.031684479828371e-01     & 3.249237036645046e-02 \\
     7.443700671611690e-01     & 4.905519963921684e-02 \\
     7.990442708271941e-01     & 5.901993373645797e-02 \\
     8.596474181980754e-01     & 6.073467932536858e-02 \\
     9.175869559770760e-01     & 5.371120494602663e-02 \\
     9.643941806993207e-01     & 3.870503119897836e-02 \\
     9.930122613589740e-01     & 1.779185041193254e-02 \\
\hline
\end{tabular}
\end{center} \end{minipage}
\begin{minipage}{0.48\linewidth} \begin{center}
\begin{tabular}{ | c | c | }
\hline
  \multicolumn{2}{|c|}{20 point quadrature rule for integrals}   \\
  \multicolumn{2}{|c|}{of the form $\int_{-1}^{1}f(x) + g(x) \log|x_{9} - x| \, dx$,} \\
  \multicolumn{2}{|c|}{where $x_{9}$ is a Gauss-Legendre node}   \\
\hline
  NODES & WEIGHTS  \\
\hline
    -9.912353127269481e-01     & 2.390220989094312e-02 \\
    -9.458275339169444e-01     & 6.859369195724087e-02 \\
    -8.534299232009863e-01     & 1.156104890379952e-01 \\
    -7.167077216635750e-01     & 1.564543949958065e-01 \\
    -5.437321547508867e-01     & 1.877085225595498e-01 \\
    -3.451233500669768e-01     & 2.074984762344433e-01 \\
    -1.328455136645940e-01     & 2.149900928447852e-01 \\
     8.079364608026202e-02     & 2.103223087093422e-01 \\
     2.840823320902124e-01     & 1.945292086962893e-01 \\
     4.667290485167077e-01     & 1.693838059093582e-01 \\
     6.204853512352519e-01     & 1.371499151597280e-01 \\
     7.394732321355052e-01     & 1.002402801599464e-01 \\
     8.201318720954396e-01     & 6.080709508468810e-02 \\
     8.607963163061316e-01     & 2.070840476545303e-02 \\
     8.727728136507039e-01     & 1.253574079839078e-02 \\
     8.916563772395616e-01     & 2.397183723558556e-02 \\
     9.192373372373420e-01     & 3.021809354380292e-02 \\
     9.500346715183706e-01     & 3.027205199814611e-02 \\
     9.775532683688947e-01     & 2.368593568061651e-02 \\
     9.954896691005256e-01     & 1.141744473788874e-02 \\
\hline
\end{tabular}
\end{center} \end{minipage}\\ \vspace{5em}
\begin{minipage}{0.48\linewidth} \begin{center}
\begin{tabular}{ | c | c | }
\hline
  \multicolumn{2}{|c|}{20 point quadrature rule for integrals}   \\
  \multicolumn{2}{|c|}{of the form $\int_{-1}^{1}f(x) + g(x) \log|x_{10} - x| \, dx$,} \\
  \multicolumn{2}{|c|}{where $x_{10}$ is a Gauss-Legendre node}   \\
\hline
  NODES & WEIGHTS  \\
\hline
    -9.857595961761246e-01     & 3.643931593123844e-02 \\
    -9.258924858821892e-01     & 8.273794370795576e-02 \\
    -8.218833662202629e-01     & 1.242680229936124e-01 \\
    -6.797470718157004e-01     & 1.586456191174843e-01 \\
    -5.075733846631832e-01     & 1.841025430283230e-01 \\
    -3.148842536644393e-01     & 1.995580161940930e-01 \\
    -1.119089520342051e-01     & 2.046841584582030e-01 \\
     9.117593460489747e-02     & 1.999093145144455e-01 \\
     2.849772954295424e-01     & 1.863566566231937e-01 \\
     4.615244999958006e-01     & 1.657233639623953e-01 \\
     6.147619568252419e-01     & 1.401105427713687e-01 \\
     7.408522006801963e-01     & 1.118164522164500e-01 \\
     8.382582352569804e-01     & 8.310260847601852e-02 \\
     9.075765988474132e-01     & 5.592493151946541e-02 \\
     9.510630103726074e-01     & 3.157199356768625e-02 \\
     9.717169387169078e-01     & 1.021414662954223e-02 \\
     9.767801773920733e-01     & 4.375212351185101e-03 \\
     9.832812993252168e-01     & 7.845621096866406e-03 \\
     9.915520723139890e-01     & 8.062764683328619e-03 \\
     9.981629455677877e-01     & 4.550772157144354e-03 \\
\hline
\end{tabular}
\end{center} \end{minipage}
\begin{minipage}{0.48\linewidth} \begin{center}
\begin{tabular}{ | c | c | }
\hline
  \multicolumn{2}{|c|}{24 point quadrature rule for integrals}   \\
  \multicolumn{2}{|c|}{of the form $\int_{0}^{1}f(x) + g(x) \log(x + \bar{x})dx$, }   \\
  \multicolumn{2}{|c|}{where $\bar{x} \ge 10^{-1} $} \\
\hline
  NODES & WEIGHTS  \\
\hline
     3.916216329415252e-02     & 4.880755296918116e-02 \\
     8.135233983530081e-02     & 3.196002785163611e-02 \\
     1.123448211344994e-01     & 3.883416642507362e-02 \\
     1.595931983965030e-01     & 5.148898992140820e-02 \\
     2.085759027831349e-01     & 4.219328148763533e-02 \\
     2.426241962027560e-01     & 3.420686213633789e-02 \\
     2.886190312538522e-01     & 5.512488680719239e-02 \\
     3.469021762354675e-01     & 6.007112809843418e-02 \\
     4.072910101569611e-01     & 6.022350479415180e-02 \\
     4.664019722595442e-01     & 5.735022004401478e-02 \\
     5.182120817844112e-01     & 4.167923417118068e-02 \\
     5.501308436771654e-01     & 3.346089628879600e-02 \\
     5.970302980854608e-01     & 5.574716218423796e-02 \\
     6.548457960388209e-01     & 5.847838243344473e-02 \\
     7.119542126106005e-01     & 5.464156990092474e-02 \\
     7.607920420946340e-01     & 4.092186343704961e-02 \\
     7.953017051155684e-01     & 3.283728166050225e-02 \\
     8.303900341517088e-01     & 3.438233273473095e-02 \\
     8.612724919009394e-01     & 3.022585192226418e-02 \\
     8.954049128027080e-01     & 3.700769701277380e-02 \\
     9.315909369155358e-01     & 3.410213679365162e-02 \\
     9.621742249068356e-01     & 2.665791885274193e-02 \\
     9.843663446380599e-01     & 1.754420526360429e-02 \\
     9.970087425823398e-01     & 7.662283104388867e-03 \\
\hline
\end{tabular}
\end{center} \end{minipage}\\
\end{center}

\begin{center}
\begin{minipage}{0.48\linewidth} \begin{center}
\begin{tabular}{ | c | c | }
\hline
  \multicolumn{2}{|c|}{24 point quadrature rule for integrals}   \\
  \multicolumn{2}{|c|}{of the form $\int_{0}^{1}f(x) + g(x) \log(x + \bar{x})dx$,}   \\
  \multicolumn{2}{|c|}{where $10^{-2} \le \bar{x} \le 10^{-1} $} \\
\hline
  NODES & WEIGHTS  \\
\hline
     1.940564616937581e-02    &  2.514022176052795e-02 \\
     4.545433992382339e-02     & 2.703526530535647e-02 \\
     7.378866604396420e-02     & 2.980872487617485e-02 \\
     1.054147718077606e-01     & 3.360626237885489e-02 \\
     1.412997888401000e-01     & 3.829678083416609e-02 \\
     1.822325567811081e-01     & 4.365651045780837e-02 \\
     2.287282121202408e-01     & 4.935846322319046e-02 \\
     2.809170925514041e-01     & 5.495967924055210e-02 \\
     3.384320962237970e-01     & 5.991162198705084e-02 \\
     4.003108031244078e-01     & 6.356960862248889e-02 \\
     4.648605571606025e-01     & 6.506868552467118e-02 \\
     5.290714994276687e-01     & 6.219588235225894e-02 \\
     5.829663557386375e-01     & 3.889986041695310e-02 \\
     6.128301889979477e-01     & 3.573431931940621e-02 \\
     6.606072156240962e-01     & 5.296315368353523e-02 \\
     7.139495966128518e-01     & 5.369033999927759e-02 \\
     7.677830914961244e-01     & 5.340793573367282e-02 \\
     8.187382423336450e-01     & 4.704756013998560e-02 \\
     8.587068551739496e-01     & 3.276576301747068e-02 \\
     8.906873285570645e-01     & 3.449175311880027e-02 \\
     9.267772492129903e-01     & 3.560168848238671e-02 \\
     9.592137652582382e-01     & 2.857367151127661e-02 \\
     9.830962712794008e-01     & 1.894042942442201e-02 \\
     9.967621546194148e-01     & 8.291994770212826e-03 \\
\hline
\end{tabular}
\end{center} \end{minipage}
\begin{minipage}{0.48\linewidth} \begin{center}
\begin{tabular}{ | c | c | }
\hline
  \multicolumn{2}{|c|}{24 point quadrature rule for integrals}   \\
  \multicolumn{2}{|c|}{of the form $\int_{0}^{1}f(x) + g(x) \log(x + \bar{x})dx$,}   \\
  \multicolumn{2}{|c|}{where $10^{-3} \le \bar{x} \le 10^{-2} $} \\
\hline
  NODES & WEIGHTS  \\
\hline
     7.571097817272427e-03     & 9.878088201321919e-03 \\
     1.800655325976786e-02     & 1.109316819462674e-02 \\
     3.003901004577040e-02     & 1.313311581321880e-02 \\
     4.462882147989575e-02     & 1.624262442061470e-02 \\
     6.295732618092606e-02     & 2.065168462990214e-02 \\
     8.644035241970913e-02     & 2.657795406825320e-02 \\
     1.166164809306920e-01     & 3.399052299072427e-02 \\
     1.546690628394902e-01     & 4.208214612865170e-02 \\
     1.999554346680615e-01     & 4.732516974042797e-02 \\
     2.434683359132119e-01     & 3.618419415803922e-02 \\
     2.800846274146029e-01     & 4.547346840583578e-02 \\
     3.368595257878888e-01     & 6.463153575242817e-02 \\
     4.044418359833648e-01     & 6.859104457897808e-02 \\
     4.685002493634456e-01     & 5.589917935916451e-02 \\
     5.185062817085154e-01     & 5.199232318335285e-02 \\
     5.811314144990846e-01     & 7.089840644422261e-02 \\
     6.545700991450585e-01     & 7.427400331494240e-02 \\
     7.276588861478224e-01     & 7.125308736931726e-02 \\
     7.960626077582168e-01     & 6.513697474660338e-02 \\
     8.572037183403355e-01     & 5.682298546820264e-02 \\
     9.091330485015775e-01     & 4.678000924507099e-02 \\
     9.503131649503738e-01     & 3.538488886617123e-02 \\
     9.795718963793163e-01     & 2.299723483013955e-02 \\
     9.961006479199827e-01     & 9.993597414733579e-03 \\
\hline
\end{tabular}
\end{center} \end{minipage}\\ \vspace{3em}
\begin{minipage}{0.48\linewidth} \begin{center}
\begin{tabular}{ | c | c | }
\hline
  \multicolumn{2}{|c|}{24 point quadrature rule for integrals}   \\
  \multicolumn{2}{|c|}{of the form $\int_{0}^{1}f(x) + g(x) \log(x + \bar{x})dx$,}   \\
  \multicolumn{2}{|c|}{where $10^{-4} \le \bar{x} \le 10^{-3} $} \\
\hline
  NODES & WEIGHTS  \\
\hline
     2.625961371586153e-03     & 3.441901737135120e-03 \\
     6.309383772392260e-03     & 3.978799794732070e-03 \\
     1.073246133489697e-02     & 4.958449505644980e-03 \\
     1.645170499644402e-02     & 6.620822501994994e-03 \\
     2.433800511777796e-02     & 9.385496468197222e-03 \\
     3.582530925992294e-02     & 1.396512052439178e-02 \\
     5.315827372101662e-02     & 2.119383832447796e-02 \\
     7.917327903614484e-02     & 3.124989308824302e-02 \\
     1.162053707416708e-01     & 4.291481168916344e-02 \\
     1.648139164451449e-01     & 5.400832278279924e-02 \\
     2.231934088488800e-01     & 6.197424674301215e-02 \\
     2.864519293820641e-01     & 6.297221626131570e-02 \\
     3.466729491189400e-01     & 5.794981636764223e-02 \\
     4.076175535528108e-01     & 6.650501614478806e-02 \\
     4.800964107543535e-01     & 7.716379373230733e-02 \\
     5.594105009204460e-01     & 8.047814122759604e-02 \\
     6.395390292352857e-01     & 7.917822434973971e-02 \\
     7.167410782176877e-01     & 7.477646096014055e-02 \\
     7.882807127957939e-01     & 6.793424765652059e-02 \\
     8.519356675821297e-01     & 5.906852968947303e-02 \\
     9.058606177202579e-01     & 4.853108558910315e-02 \\
     9.485539755760567e-01     & 3.666228059710319e-02 \\
     9.788566874094059e-01     & 2.380850649522536e-02 \\
     9.959649506960162e-01     & 1.034186239262945e-02 \\
\hline
\end{tabular}
\end{center} \end{minipage}
\begin{minipage}{0.48\linewidth} \begin{center}
\begin{tabular}{ | c | c | }
\hline
  \multicolumn{2}{|c|}{24 point quadrature rule for integrals}   \\
  \multicolumn{2}{|c|}{of the form $\int_{0}^{1}f(x) + g(x) \log(x + \bar{x})dx$,}   \\
  \multicolumn{2}{|c|}{where $10^{-5} \le \bar{x} \le 10^{-4} $} \\
\hline
  NODES & WEIGHTS  \\
\hline
     7.759451679242260e-04     & 1.049591733965263e-03 \\
     1.952854410117286e-03     & 1.314968855711329e-03 \\
     3.429053832116395e-03     & 1.651475072547296e-03 \\
     5.301128540262913e-03     & 2.135645684467029e-03 \\
     7.878118775220067e-03     & 3.165043382856636e-03 \\
     1.205537050949829e-02     & 5.479528688655274e-03 \\
     1.965871512055557e-02     & 1.028817002915096e-02 \\
     3.403328641997047e-02     & 1.923291785614007e-02 \\
     5.947430305925957e-02     & 3.212643438782854e-02 \\
     9.873500543531440e-02     & 4.638626850049229e-02 \\
     1.518862681939413e-01     & 5.960676923068444e-02 \\
     2.171724325134259e-01     & 7.052360405410943e-02 \\
     2.919941878735093e-01     & 7.863451090237836e-02 \\
     3.734637353255530e-01     & 8.381771698595157e-02 \\
     4.586710018443288e-01     & 8.612755554083525e-02 \\
     5.448057416999684e-01     & 8.569938467103264e-02 \\
     6.292158981939618e-01     & 8.271051499695768e-02 \\
     7.094415843889587e-01     & 7.736692567834522e-02 \\
     7.832417328632321e-01     & 6.990012937760461e-02 \\
     8.486194141302759e-01     & 6.056687669667680e-02 \\
     9.038469149367938e-01     & 4.964868706783169e-02 \\
     9.474898150194623e-01     & 3.745026957972177e-02 \\
     9.784290662963747e-01     & 2.429741981889855e-02 \\
     9.958843370550371e-01     & 1.054906616108520e-02 \\
\hline
\end{tabular}
\end{center} \end{minipage}\\
\end{center}

\begin{center}
\begin{minipage}{0.48\linewidth} \begin{center}
\begin{tabular}{ | c | c | }
\hline
  \multicolumn{2}{|c|}{24 point quadrature rule for integrals}   \\
  \multicolumn{2}{|c|}{of the form $\int_{0}^{1}f(x) + g(x) \log(x + \bar{x})dx$,}   \\
  \multicolumn{2}{|c|}{where $10^{-6} \le \bar{x} \le 10^{-5} $} \\
\hline
  NODES & WEIGHTS  \\
\hline
     3.126377187332637e-04     & 4.136479682893960e-04 \\
     7.671264269072188e-04     & 5.068714387414649e-04 \\
     1.359575160544077e-03     & 7.008932527842778e-04 \\
     2.238313285727558e-03     & 1.110264922990352e-03 \\
     3.770276623583326e-03     & 2.120108385941761e-03 \\
     7.146583956092048e-03     & 5.249076343206215e-03 \\
     1.635515250548719e-02     & 1.450809938905405e-02 \\
     3.828062855101241e-02     & 2.987724029376343e-02 \\
     7.628984500206759e-02     & 4.593298717863718e-02 \\
     1.294255336121595e-01     & 5.987634475538021e-02 \\
     1.949876755761554e-01     & 7.065953519392547e-02 \\
     2.693852297828856e-01     & 7.729918562776261e-02 \\
     3.469762441631538e-01     & 7.556635340171830e-02 \\
     4.122748928895491e-01     & 5.234123638339037e-02 \\
     4.662499202239145e-01     & 6.532130125393047e-02 \\
     5.421402737123784e-01     & 8.188272080198840e-02 \\
     6.248832413655412e-01     & 8.237354882288161e-02 \\
     7.053258496784840e-01     & 7.795795664563893e-02 \\
     7.798841313231049e-01     & 7.076514272025076e-02 \\
     8.461534275163378e-01     & 6.145788741452406e-02 \\
     9.022312524979976e-01     & 5.044339641339403e-02 \\
     9.465899812310277e-01     & 3.807817118430632e-02 \\
     9.780549563823810e-01     & 2.471549011101626e-02 \\
     9.958125149101927e-01     & 1.073289672726758e-02 \\
\hline
\end{tabular}
\end{center} \end{minipage}
\begin{minipage}{0.48\linewidth} \begin{center}
\begin{tabular}{ | c | c | }
\hline
  \multicolumn{2}{|c|}{24 point quadrature rule for integrals}   \\
  \multicolumn{2}{|c|}{of the form $\int_{0}^{1}f(x) + g(x) \log(x + \bar{x})dx$,}   \\
  \multicolumn{2}{|c|}{where $10^{-7} \le \bar{x} \le 10^{-6} $} \\
\hline
  NODES & WEIGHTS  \\
\hline
     1.019234906342863e-04     & 1.349775051746596e-04 \\
     2.506087227631447e-04     & 1.663411550150506e-04 \\
     4.461429005344285e-04     & 2.328782111562424e-04 \\
     7.422845421202523e-04     & 3.804721779784063e-04 \\
     1.289196091156456e-03     & 7.930350452911450e-04 \\
     2.739287668024851e-03     & 2.600694722423854e-03 \\
     9.075168969969708e-03     & 1.212249113599252e-02 \\
     2.968005234555358e-02     & 2.946708975720586e-02 \\
     6.781742979962609e-02     & 4.647771960691390e-02 \\
     1.217792474402805e-01     & 6.095376889009233e-02 \\
     1.886625378438471e-01     & 7.224844725827559e-02 \\
     2.650602155844836e-01     & 7.986429603884565e-02 \\
     3.465113608339080e-01     & 8.143206462900546e-02 \\
     4.178374197420536e-01     & 5.040529357007135e-02 \\
     4.597624982511183e-01     & 5.592137651001418e-02 \\
     5.348065111487157e-01     & 8.398073572656715e-02 \\
     6.194640153146728e-01     & 8.402586870225486e-02 \\
     7.013481004172354e-01     & 7.922223490159952e-02 \\
     7.770386175609082e-01     & 7.177919251691964e-02 \\
     8.442211768916794e-01     & 6.227551999401272e-02 \\
     9.010272836291835e-01     & 5.108407212719758e-02 \\
     9.459409782755001e-01     & 3.854783279333592e-02 \\
     9.777905486554876e-01     & 2.501496650831813e-02 \\
     9.957622871041650e-01     & 1.086176801402067e-02 \\
\hline
\end{tabular}
\end{center} \end{minipage}\\ \vspace{3em}
\begin{minipage}{0.48\linewidth} \begin{center}
\begin{tabular}{ | c | c | }
\hline
  \multicolumn{2}{|c|}{24 point quadrature rule for integrals}   \\
  \multicolumn{2}{|c|}{of the form $\int_{0}^{1}f(x) + g(x) \log(x + \bar{x})dx$,}   \\
  \multicolumn{2}{|c|}{where $10^{-8} \le \bar{x} \le 10^{-7} $} \\
\hline
  NODES & WEIGHTS  \\
\hline
     3.421721832247593e-05     & 4.559730842497453e-05 \\
     8.533906255442380e-05     & 5.840391255974745e-05 \\
     1.563524616155011e-04     & 8.761580900682040e-05 \\
     2.746612401575526e-04     & 1.617264666294872e-04 \\
     5.408643931265062e-04     & 4.433543035169213e-04 \\
     1.782382096488333e-03     & 3.116175111368442e-03 \\
     1.101243912052365e-02     & 1.655494413772595e-02 \\
     3.553172024884285e-02     & 3.242539256461602e-02 \\
     7.554170435463801e-02     & 4.734426463929677e-02 \\
     1.295711894941649e-01     & 6.032614603579952e-02 \\
     1.953213037793089e-01     & 7.069975187373848e-02 \\
     2.699680545714222e-01     & 7.806973621204365e-02 \\
     3.503697281371090e-01     & 8.216350598137868e-02 \\
     4.330838596494367e-01     & 8.261286657092808e-02 \\
     5.141801680435878e-01     & 7.883476216668445e-02 \\
     5.895097016206093e-01     & 7.157205125318401e-02 \\
     6.582708672338614e-01     & 6.703064468754417e-02 \\
     7.252543617887320e-01     & 6.706137273719630e-02 \\
     7.914154485613720e-01     & 6.449984116349734e-02 \\
     8.528383935857844e-01     & 5.775434959088197e-02 \\
     9.059696536862878e-01     & 4.812600239023880e-02 \\
     9.484664124578303e-01     & 3.661415869304224e-02 \\
     9.787863313133854e-01     & 2.386304203446463e-02 \\
     9.959482975155097e-01     & 1.038268695581411e-02 \\
\hline
\end{tabular}
\end{center} \end{minipage}
\begin{minipage}{0.48\linewidth} \begin{center}
\begin{tabular}{ | c | c | }
\hline
  \multicolumn{2}{|c|}{24 point quadrature rule for integrals}   \\
  \multicolumn{2}{|c|}{of the form $\int_{0}^{1}f(x) + g(x) \log(x + \bar{x})dx$,}   \\
  \multicolumn{2}{|c|}{where $10^{-9} \le \bar{x} \le 10^{-8} $} \\
\hline
  NODES & WEIGHTS  \\
\hline
     6.538987938840374e-06     & 1.500332421093607e-05 \\
     2.613485075847413e-05     & 2.367234654253158e-05 \\
     5.664183720634991e-05     & 4.007286246706405e-05 \\
     1.179374114362569e-04     & 9.497743501485505e-05 \\
     3.299119431334128e-04     & 4.619067037944727e-04 \\
     3.626828607577001e-03     & 9.985382463808036e-03 \\
     2.265102906572155e-02     & 2.805741744607257e-02 \\
     5.896796231680340e-02     & 4.404106103008398e-02 \\
     1.092496277855923e-01     & 5.548413172821072e-02 \\
     1.666701689499393e-01     & 5.693235996372726e-02 \\
     2.196889385898800e-01     & 5.087307376046002e-02 \\
     2.770352260035617e-01     & 6.593729718379782e-02 \\
     3.483163928268329e-01     & 7.335680008972614e-02 \\
     4.153287664837260e-01     & 5.675029500743735e-02 \\
     4.695624219668608e-01     & 6.117926027541254e-02 \\
     5.421129318998841e-01     & 8.004805067067550e-02 \\
     6.238832212055707e-01     & 8.196991767042605e-02 \\
     7.041842972237081e-01     & 7.800219127200407e-02 \\
     7.788817007552110e-01     & 7.097175077519494e-02 \\
     8.453877637047045e-01     & 6.171193295041172e-02 \\
     9.017178251963006e-01     & 5.068671319716005e-02 \\
     9.462999385952402e-01     & 3.827738423897266e-02 \\
     9.779333485180249e-01     & 2.485063762733620e-02 \\
     9.957890687155009e-01     & 1.079284973329516e-02 \\
\hline
\end{tabular}
\end{center} \end{minipage}\\
\end{center}

\begin{center}
\begin{minipage}{0.48\linewidth} \begin{center}
\begin{tabular}{ | c | c | }
\hline
  \multicolumn{2}{|c|}{24 point quadrature rule for integrals}   \\
  \multicolumn{2}{|c|}{of the form $\int_{0}^{1}f(x) + g(x) \log(x + \bar{x})dx$,}   \\
  \multicolumn{2}{|c|}{where $10^{-10} \le \bar{x} \le 10^{-9} $} \\
\hline
  NODES & WEIGHTS  \\
\hline
     6.725520559705825e-06     & 8.128391913974039e-05 \\
     6.986424152770461e-06    & -7.773900735768282e-05 \\
     1.217363416714366e-05     & 1.287386499666193e-05 \\
     2.677746219601529e-05     & 1.895577251914526e-05 \\
     5.597036348896741e-05     & 4.732580352158076e-05 \\
     2.729343280943077e-04     & 9.857909615386162e-04 \\
     9.445526806263141e-03     & 1.756872897270054e-02 \\
     3.556725025161542e-02     & 3.439422017906772e-02 \\
     7.765556668177810e-02     & 4.944188361792970e-02 \\
     1.336848150648662e-01     & 6.219733934997792e-02 \\
     2.011576917683550e-01     & 7.228007436918939e-02 \\
     2.772736854314979e-01     & 7.944986391225688e-02 \\
     3.590124362607926e-01     & 8.347646288178011e-02 \\
     4.430074035214462e-01     & 8.380433020121207e-02 \\
     5.247388219574510e-01     & 7.832768209682506e-02 \\
     5.961053238782420e-01     & 6.300796225242940e-02 \\
     6.547331131213409e-01     & 5.923406014585053e-02 \\
     7.192258519628951e-01     & 6.834293563803810e-02 \\
     7.874251789073102e-01     & 6.660337204499726e-02 \\
     8.505852012775045e-01     & 5.911988751082552e-02 \\
     9.047824617894323e-01     & 4.893575310568894e-02 \\
     9.479045131744448e-01     & 3.708256438629509e-02 \\
     9.785770588866582e-01     & 2.411463784693618e-02 \\
     9.959104692340199e-01     & 1.048087156697020e-02 \\
\hline
\end{tabular}
\end{center} \end{minipage}
\begin{minipage}{0.48\linewidth} \begin{center}
\begin{tabular}{ | c | c | }
\hline
  \multicolumn{2}{|c|}{24 point quadrature rule for integrals}   \\
  \multicolumn{2}{|c|}{of the form $\int_{0}^{1}f(x) + g(x) \log(x + \bar{x})dx$,}   \\
  \multicolumn{2}{|c|}{where $10^{-11} \le \bar{x} \le 10^{-10} $} \\
\hline
  NODES & WEIGHTS  \\
\hline
     2.828736694877886e-08     & 1.665602686704325e-05 \\
     2.302233157554212e-06     & 2.577419924039251e-06 \\
     5.853587143444178e-06     & 4.957941112780975e-06 \\
     1.451588770083244e-05     & 1.537074702915107e-05 \\
     9.711965099273031e-05     & 4.640075239797995e-04 \\
     9.004761967373848e-03     & 1.705687938176189e-02 \\
     3.442077924035546e-02     & 3.349724914160473e-02 \\
     7.543926781582543e-02     & 4.820210872119093e-02 \\
     1.300373356318913e-01     & 6.054547286337976e-02 \\
     1.955182772803384e-01     & 6.984354388121057e-02 \\
     2.683608546664295e-01     & 7.498721497014774e-02 \\
     3.430029178740901e-01     & 7.240620145057083e-02 \\
     4.085056107803621e-01     & 5.774925310174693e-02 \\
     4.660198270439085e-01     & 6.238505554837956e-02 \\
     5.336124745634699e-01     & 6.940394677081842e-02 \\
     5.985245800106473e-01     & 5.910843483407385e-02 \\
     6.564089719608276e-01     & 6.059752321454190e-02 \\
     7.216666024232565e-01     & 6.823362237770209e-02 \\
     7.893712241343741e-01     & 6.593839664071163e-02 \\
     8.518883782001418e-01     & 5.853014420243146e-02 \\
     9.055688088881344e-01     & 4.849217100974983e-02 \\
     9.483163097840529e-01     & 3.677417821170115e-02 \\
     9.787413692715607e-01     & 2.392585642844202e-02 \\
     9.959413203611228e-01     & 1.040149939671874e-02 \\
\hline
\end{tabular}
\end{center} \end{minipage}\\ \vspace{2em}
\begin{minipage}{0.48\linewidth} \begin{center}
\begin{tabular}{ | c | c | }
\hline
  \multicolumn{2}{|c|}{24 point quadrature rule for integrals}   \\
  \multicolumn{2}{|c|}{of the form $\int_{0}^{1}f(x) + g(x) \log(x + \bar{x})dx$,}   \\
  \multicolumn{2}{|c|}{where $10^{-12} \le \bar{x} \le 10^{-11} $} \\
\hline
  NODES & WEIGHTS  \\
\hline
     6.147063879573664e-07     & 8.763741095000331e-07 \\
     2.102921984985835e-06     & 1.784696796288373e-05 \\
     2.188366117432289e-06    & -1.795398395983826e-05 \\
     3.482602942694880e-06     & 5.117514567175025e-06 \\
     2.768001888608636e-05     & 1.698863549284390e-04 \\
     8.942779215792784e-03     & 1.701975216672032e-02 \\
     3.432218364237253e-02     & 3.346025972593909e-02 \\
     7.530931328026620e-02     & 4.817949622196712e-02 \\
     1.298983048592572e-01     & 6.055152664710045e-02 \\
     1.954020797117703e-01     & 6.988313730886592e-02 \\
     2.682970870436427e-01     & 7.504602275463067e-02 \\
     3.429540704041702e-01     & 7.230942674874111e-02 \\
     4.080399755202422e-01     & 5.705952259766429e-02 \\
     4.652562798154792e-01     & 6.265021180818162e-02 \\
     5.333220999210325e-01     & 6.993669694523695e-02 \\
     5.986982369433125e-01     & 5.937130986945129e-02 \\
     6.564773600603511e-01     & 6.026572020863567e-02 \\
     7.215159032030418e-01     & 6.815292696374753e-02 \\
     7.892098210760941e-01     & 6.596804590657802e-02 \\
     8.517672777806986e-01     & 5.857483758149194e-02 \\
     9.054906995605498e-01     & 4.853209199396977e-02 \\
     9.482736017320823e-01     & 3.680469214176019e-02 \\
     9.787238593479314e-01     & 2.394561701705853e-02 \\
     9.959379852805677e-01     & 1.041005152890511e-02 \\
\hline
\end{tabular}
\end{center} \end{minipage}
\begin{minipage}{0.48\linewidth} \begin{center}
\begin{tabular}{ | c | c | }
\hline
  \multicolumn{2}{|c|}{24 point quadrature rule for integrals}   \\
  \multicolumn{2}{|c|}{of the form $\int_{0}^{1}f(x) + g(x) \log(x + \bar{x})dx$,}   \\
  \multicolumn{2}{|c|}{where $10^{-13} \le \bar{x} \le 10^{-12} $} \\
\hline
  NODES & WEIGHTS  \\
\hline
     4.523740015216508e-08     & 4.418138082366788e-07 \\
     4.281855233588279e-07     & 4.389108058643120e-07 \\
     1.036900153156159e-06     & 9.539585150737866e-07 \\
     7.825849325746907e-06     & 5.823980947200484e-05 \\
     8.617419723953112e-03     & 1.634464263521301e-02 \\
     3.268881163637599e-02     & 3.129682188728318e-02 \\
     6.988441391437043e-02     & 4.212468617589480e-02 \\
     1.142202307676442e-01     & 4.505120897719191e-02 \\
     1.596471081833281e-01     & 4.769069780026684e-02 \\
     2.135336418959620e-01     & 6.038503382768951e-02 \\
     2.781100275296151e-01     & 6.695343672694180e-02 \\
     3.433392803364457e-01     & 6.163298712826237e-02 \\
     4.019960595528027e-01     & 5.877742624357513e-02 \\
     4.656415679416787e-01     & 6.800053637773440e-02 \\
     5.334880548894250e-01     & 6.516918103589647e-02 \\
     5.943298528903542e-01     & 5.853785375926075e-02 \\
     6.562968737815924e-01     & 6.639396325654251e-02 \\
     7.250343344601498e-01     & 6.948738324081696e-02 \\
     7.928820737781136e-01     & 6.538801703374268e-02 \\
     8.546103048745466e-01     & 5.761503751629250e-02 \\
     9.073762310762705e-01     & 4.761344859555310e-02 \\
     9.493253659835347e-01     & 3.607033097268266e-02 \\
     9.791606801267259e-01     & 2.345690720840071e-02 \\
     9.960217573957566e-01     & 1.019557402722854e-02 \\
\hline
\end{tabular}
\end{center} \end{minipage}\\
\end{center}

\begin{center}
\begin{tabular}{ | c | c | }
\hline
  \multicolumn{2}{|c|}{24 point quadrature rule for integrals}   \\
  \multicolumn{2}{|c|}{of the form $\int_{0}^{1}f(x) + g(x) \log(x + \bar{x})dx$,}   \\
  \multicolumn{2}{|c|}{where $10^{-14} \le \bar{x} \le 10^{-13} $} \\
\hline
  NODES & WEIGHTS  \\
\hline
     6.025980282801020e-08     & 9.079353616441234e-07 \\
     6.411245262925473e-08    & -8.390389042773805e-07 \\
     1.862815529429129e-07     & 2.782460677485016e-07 \\
     2.029190208906422e-06     & 1.821115881362725e-05 \\
     8.902881307076499e-03     & 1.695809650660321e-02 \\
     3.420089035164912e-02     & 3.336370146025145e-02 \\
     7.508687525931594e-02     & 4.807898681796971e-02 \\
     1.295858123029775e-01     & 6.047672723211479e-02 \\
     1.950409815188335e-01     & 6.986774906175534e-02 \\
     2.679751967812604e-01     & 7.515608233194288e-02 \\
     3.428525062164689e-01     & 7.264249904037610e-02 \\
     4.080941369413548e-01     & 5.672507168477261e-02 \\
     4.646644511900009e-01     & 6.220316364524964e-02 \\
     5.328071517215501e-01     & 7.032362652293805e-02 \\
     5.978508749698001e-01     & 5.742730804758014e-02 \\
     6.521214523350964e-01     & 5.644075454541152e-02 \\
     7.134921670665336e-01     & 6.318643666150391e-02 \\
     7.679317896479284e-01     & 3.945995610428228e-02 \\
     8.029718487208403e-01     & 4.324200884758527e-02 \\
     8.551101435866935e-01     & 5.478223695609097e-02 \\
     9.067319102017767e-01     & 4.740856250832772e-02 \\
     9.487765213293372e-01     & 3.633314063504751e-02 \\
     9.788979796532736e-01     & 2.372788917088821e-02 \\
     9.959684838634199e-01     & 1.033036588606145e-02 \\
\hline
\end{tabular}
\end{center}

\lsp

\end{document}